\documentclass[preprint,10pt,authoryear]{elsarticle}
\usepackage{natbib}
\setcitestyle{square,sort,comma,numbers}
\usepackage{amsmath}
\usepackage{graphicx}
\usepackage{color}
\usepackage{cases}
\usepackage{float}
\usepackage{amssymb}
\usepackage{subcaption}
\usepackage{multirow}
\usepackage{lscape}
\usepackage[utf8x]{inputenc}
\usepackage[T1]{fontenc}
\usepackage{textcomp} 
\usepackage{graphicx}
\usepackage{ucs}
\usepackage{amssymb}
\usepackage{amsmath}
\usepackage{amsfonts}
\usepackage{latexsym}
\usepackage{hyperref}
\hypersetup{colorlinks=true}
\hypersetup{urlcolor=blue}
\hypersetup{linkcolor=red}
\hypersetup{citecolor=blue}
\usepackage{upgreek}
\usepackage{graphics}
\usepackage{float}
\usepackage{caption}
\usepackage{epsfig}
\usepackage{epstopdf}
\usepackage{caption}
\usepackage{subcaption}
\usepackage{psfrag}
\usepackage{titlesec}
\newcommand{\sbt}{\,\begin{picture}(-1,1)(-1,-3)\circle*{3}\end{picture}\;\,}

\newcommand{\vett}[1]{\boldsymbol{#1}}

\titleformat{\paragraph}
{\normalfont\normalsize\bfseries}{\theparagraph}{1em}{}
\titlespacing*{\paragraph}
{0pt}{3.25ex plus 1ex minus .2ex}{1.5ex plus .2ex}

\begin{document}
\begin{frontmatter}
\small
\title{A hybrid MGA-MSGD ANN training approach for approximate solution of linear elliptic PDEs}  

\author[unilux]{\corref{cor2}H. Dehghani} 
\ead{hamidreza.dehghani@uni.lu}
\author[unilux]{A. Zilian}
\cortext[cor2]{Corresponding author}
\address[unilux]{Institute of Computational Engineering and Sciences, Faculty of Science, Technology and Medicine, University of Luxembourg, 6 Avenue de la Fonte, 4364 Esch-sur-Alzette, Luxembourg.}

\begin{abstract}

We introduce a hybrid "Modified Genetic Algorithm-Multilevel Stochastic Gradient Descent" (MGA-MSGD) training algorithm that considerably improves accuracy and efficiency of solving 3D mechanical problems described, in strong-form, by PDEs via ANNs (Artificial Neural Networks). This presented approach allows the selection of a number of \textit{locations of interest} at which the state variables are expected to fulfil the governing equations associated with a physical problem. Unlike classical PDE approximation methods  such  as  finite differences or the finite element method,  there  is no need to establish and reconstruct the physical field quantity throughout the computational domain in order to predict the mechanical response  at  specific  locations  of  interest. The basic idea of MGA-MSGD is the manipulation of the learnable parameters' components responsible for the error explosion so that we can train the network with relatively larger learning rates which avoids trapping in local minima. The proposed training approach is less sensitive to the learning rate value, training points density and distribution, and the random initial parameters. The distance function to minimise is where we introduce the PDEs including any physical laws and conditions (so-called, Physics Informed ANN). The Genetic algorithm is modified to be suitable for this type of ANN in which a Coarse-level Stochastic Gradient Descent (CSGD) is exploited to make the decision of the offspring qualification. Employing the presented approach, a considerable improvement in both accuracy and efficiency, compared with standard training algorithms such classical SGD and Adam optimiser, is observed. The local displacement accuracy is studied and ensured by introducing the results of Finite Element Method (FEM) at sufficiently fine mesh as the reference displacements. A slightly more complex problem is solved ensuring the feasibility of the methodology.
\end{abstract}

\end{frontmatter}

\paragraph{Declaration of interest}
\textit{Declarations of interest: none}
\section{Introduction}
The growth in computing power together with data science has allowed the emergence of data-driven computational mechanics in which one aims at e.g. material parameter identification, improving or complementing standard computational mechanics approaches such as FEM, or governing equations discovery, etc. \citep{HDAZ2020, OISHI2017327, KIRCHDOERFER201681, Raissi2018}.
Over the past few years, Artificial Neural Networks (ANNs) \citep{Rosenblatt58theperceptron},  which provide input-output transfer means for complex problems, have been applied successfully in this field. One major issue of using this approach is the training procedure which can be considerably time-consuming and sensitive to several ANN-related parameters. In the latter process, the learnable parameters (the ones such as its \emph{weights} and \emph{biases} that determine the ANN output) are tuned by minimising a distance function called the cost function which is calculated based on, in the simplest case, a dataset including the exact solutions to be interpolated. 
Furthermore, ANN training sometimes is only a one-time calculation \citep{HDAZ2020} which allows the choice of non-efficient approaches. 

However, for \emph{strong solution} of mechanical problems (i.e. the solution of the strong-form of the problem) described by a system of Partial Differential Equations (PDEs) a training dataset of exact solutions is not available in which case the cost function is constructed based on the governing PDEs (so-called Physics Informed ANNs) \citep{Raissi2018}. In fact, in this concept, solving a problem such as a linear elastic  deformation problem by ANN means training a network which provides a space-displacement relationship satisfying the strong form of governing physical laws in the form of PDEs which is enforced via the cost function. For a better understanding of the general methodology one could assume that ANN is a complex ansatz, relating the coordinates of interest to the corresponding displacements, whose parameters are tuned by minimising the mentioned distance function.
The provided solution (displacement field), as well as its spatial gradient (strain/stress field), are continuous, directly differentiable, and free of space discretisation. The latter, results in increased flexibility and less complexity of the problem implementation. For example, for obtaining the mechanical response at a specific coordinate one does not need to discretise and solve the entire domain, which will be shown in the results. However, the training procedure can be challenging and inefficient as we need to compute the third order gradients to solve every problem entirely by ANN without having a training dataset that provides the exact results.

There are several training methods derived from Back Propagation (BP) such as classical Stochastic Gradient Descent (SGD), Adam, and Adagrad \citep{Robbins2007ASA,kingma2014adam, Duchi2011} which train the network based on a graph of the derivatives of the cost function with respect to each learnable parameter. Besides, in the present problem, the cost function is obtained by calculating the second order partial derivatives of every output with respect to every input via chain rule which implies calculation of the  third order derivatives for the Back Propagation (BP) process. Having a large number of training points and intermediate/ANN parameters, SGD does not provide acceptable efficiency and, most of the time, it does not provide an acceptable outcome for 3D problems.
Moreover, it is shown that using the SGD-based approaches, the high number of iterations necessary to reach an acceptable accuracy is a major issue in 1D and 2D problems \citep{Dockhorn2019}. This issue imposes a limit on the choice of the size of training inputs and the network's parameters which determines the accuracy and robustness of the solution if we make use of the mentioned training approaches. Another important issue is that the SGD based algorithms often get trapped in local minima as a consequence of the lack of global search. Moreover, they are sensitive to the initial choice of the network's parameters which are, normally, chosen randomly. On the other hand, in this application, these methods have a better potential of the local search compared to Evolutionary Algorithms (EAs) such as Genetic Algorithm (GA).

EAs are population-based optimisation and search algorithms inspired by natural/biological evolution mechanisms such as survival of the fittest \citep{Yao1999, Vikhar2016}. These methods have been widely employed for global optimisation problems where the connectionist approaches such as SGD are not effective. EAs include several evolution strategies such as Evolutionary Programming (EP) and Genetic Algorithm (GA) \citep{Fogel1998, Yao1999, Holland1992, Goldberg1989}.  These methods are more likely to avoid local minima because of the strategies such as \emph{mutation} that are responsible for carrying a global search and maintaining the diversity of the population. Traditionally, they are not gradient-based methods although, due to the type of the presented problem, the gradient plays a role here (to compute the loss). Moreover, EAs are less sensitive to the initial choice of the parameters. However, in the case of ANN's weights and biases training these approaches always search for a globally optimal solution and are usually inefficient in finding fine-tuned settings. In practice, in the case of ANN training, EAs are efficient in finding good initial parameters but a local search algorithm is needed to find the optimal setting.

Hybrid training algorithms have been preferred and successfully applied in several scenarios of interest from civil engineering and material science to pattern, strategy, image, and text recognition \citep{Adeli1994, Skinner1995, taha1995evolutionary, Lee1996,Aguilar}. In these methods, usually, GA is used with the goal to find a set of initial learnable parameters that is in the basin of attraction of the global minimum and a connectionist search algorithm such as SGD is employed to take fine steps to \emph{go down the valley} and minimise the distance function. In other words, GA identifies the deepest \emph{valley} and SGD determines a setting at or close to the minimum. In practice, in several studies, the hybrid algorithms are reported to outperform either GA or BP alone \citep{Belew90evolvingnetworks, TOPCHY1997240, yan1997}. The first part of the hybrid methods (GA) is based on choosing several (usually hundreds) random initial settings of learnable parameters (each one is called an individual which together form a population), compute the error/cost function for each one, affix a score called \emph{fitness} based on their distance from the target, select a certain number of individuals (surviving individuals), and construct a new population that inherits features from the latter. This procedure is repeated until a population with acceptable fitnesses is achieved, which will be explained in details. As an output, an initial setting of learnable parameters is obtained, which has a smaller distance from the target which explains the better efficiency of the algorithm. However, we show that starting from the latter initial point in the space of learnable parameters does not necessarily lead us to a smaller final distance function. Moreover, in the present problem, as we need to compute the second gradients to evaluate an individual, it is a time-consuming procedure to calculate the fitnesses for a population in each iteration. Consequently, there is a need to introduce a new training method for an accurate and efficient problem-solving procedure based on a small network with a small number of training points and iterations which is able to solve problems in 3D.

A novel hybrid method called Modified Genetic Algorithm-Multilevel Stochastic Gradient Descent (MGA-MSGD) is introduced in the present study which is, for the first time, able to provide the continuous solution to the strong form of 3D PDEs in solid mechanics (here, a linear elastic problem) efficiently and accurately. In this framework, the parameters that cause the "error explosion" (sharp unbound increase in cost function) in the training procedure are constantly modified so that a large learning rate can be used which helps to avoid several local minima. In the MGA part, the standard GA is modified in which, for example, the individuals are the binary representation of each learnable parameters and the population consists of all the weights and biases so that a direct manipulation of the parameters is done. We introduce a new way of selection based on the computation of the \emph{importance} of the individuals. In this method, for better efficiency, we only consider the qualified offsprings to produce a new generation so that every generation is necessarily \emph{better} than the previous one. The population qualification procedure is based on a coarse-level SGD training which indicates if we are at a basin of attraction of a global (or near global) minimum. We show that this approach is considerably more reliable than the conventional evaluation of the individuals and population. At the end of MGA-CSGD (Coarse-level SGD), we introduce a fine-scale SGD to step closer to the minimum with considerably smaller steps which completes the procedure of MGA-MSGD method. This approach is less sensitive to the initial parameters, the number and distribution of the training points, and learning rates. We compare the results of MGA-MSGD with the ones of a simple SGD training to demonstrate its effectiveness and we measure its accuracy using the approximate response provided by the Finite Element Method (FEM). 

In the following section, the ANN setup including its architecture, activation and cost functions, and the implementation of the physical laws (in the form of PDEs) and Dirichlet and Neumann Boundary Conditions (BCs) etc. is explained. The details of the novel training approach MGA-MSGD are presented in Section \ref{sec_MGAMSGD}. A thorough sensitivity analysis follows this in Section \ref{sec_sensitivity}, which provides a deep understanding of each parameter, its effects on the results, and the advantages of the proposed training method. Finally, in Section \ref{sec_conclusion} the concluding remarks are provided.

\section{ANN and continuum mechanics} \label{sec_ANNsetup}
\begin{figure}
	\includegraphics[width=12cm]{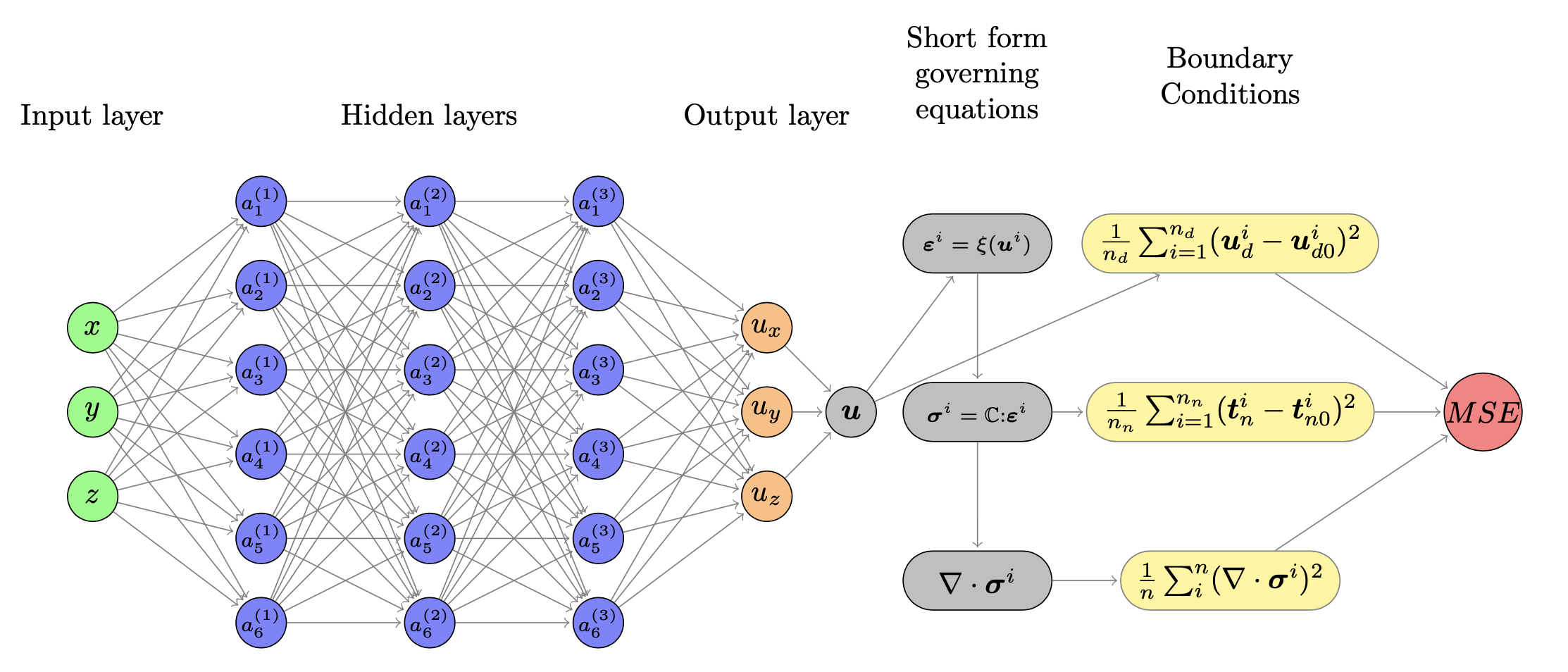}
	\centering
	\caption{Schematic representation of ANN setup for solving PDEs in solid mechanics with three spatial dimensions. Note that $\vett u_d$ represents the displacement vectors of the points located on $\Gamma_d$ and $\vett{t}_n$ indicates the Cauchy stress vector of the points on $\Gamma_n$.}
	\label{fig_ANN}
\end{figure}

In this section, we provide details of the ANN for continuous solution of a continuum mechanics problem which is assumed to be a 3D linear elastic solid deformation problem in a displacement-based formulation. 
The inputs of the ANN are the spatial coordinates  $(x,y,z)$ of every point of interest inside the domain $\Omega \subset \mathbb{R}^3$ and the outputs are their corresponding components of the displacement vector  $\vett u =(u_x,u_y,u_z)$. Such a setting is shown in Figure \ref{fig_ANN} and can be written in the form of the following transfer function
\begin{equation}
f{:}( x, y, z) \rightarrow u_x, u_y, u_z
\end{equation}

Here, we make use of a network with $N_h$ hidden layers and $N_{nh}$ neurones in each layer. The following activation function
\begin{equation}
 ELU(\sbt) = \max (0,\sbt)+ \min (0, \alpha * (\exp (\sbt) - 1))
 \end{equation}
 with the coefficient $\alpha = 1$ as we need the 2nd order derivative of the outputs (displacements) with respect to the inputs (coordinates) as well as the 3rd order derivative with respect to the learnable parameters for the ANN training procedure.
The latter starts from, in this case, a linear Feed Forward expressed as 
  \begin{align} 
 \psi^{(i)}_k&=w^{(i)}_{kj}a^{(i-1)}_j+b^{(i)}_k \label{ANN1} \\
 a^{(i)}_k&=ELU(\psi^{(i)}_k) \label{activeF}
 \end{align}
 where $a^{(i)}_k$ is the value of the $k$ -th neurone in the $(i)$ -th hidden layer with $a^{(0)}_k$ being the coordinates of interest ($\vett a^{(0)} = (x,y,z)$). $ i \in(1,..,N_h)$, $k$ and $j$ indicate the number of the neurones in the $(i)$ -th and $(i-1)$ -th layers, respectively. The ANN's outputs are then calculated by
 \begin{equation}
 u_k = w^{(N_h+1)}_{kj} a^{(N_h)}_j. \label{ANN3}
 \end{equation}

Equations \eqref{ANN1}-\eqref{ANN3} provide a complex ansatz relating spatial coordinates to their corresponding displacements.

Having a prediction of the displacements ($u_k$), we need to introduce an appropriate distance function to quantify the deviation of the solution provided by ANN from the exact/unknown displacements satisfying the following continuum mechanics governing equations.

\paragraph{Residuals of the strong-form governing equations in the domain}
The mechanical deformation problem considered here is based on the classical continuum mechanics of linear elasticity. The differential form of the governing equations are presented here as residual expression assuming the presence of approximate solutions to the unknown physical fields: displacement vector $\vett u(\vett x)$, strain tensor $\vett{\varepsilon}(\vett x)$, and stress tensor $\vett \sigma(\vett x)$ in the domain, and displacement vector $\vett{u}(x)$ and stress vector $\vett t(x)$ on the boundary.

The residual of the linearised strain tensor (assuming infinitesimal strain)
\begin{equation}
\vett r_1(\vett x) = \vett{\varepsilon}(\vett x) - \frac{1}{2} \left(\nabla \vett u(\vett x) + \nabla \vett u^T(\vett x)\right) \quad \quad \textrm{in} \quad \quad \Omega
\end{equation}
is a function of the (approximate) displacement state $\vett u(\vett x)$.
The residual of the linear momentum equation, that describes the equilibrium, is 
\begin{equation}
\vett r_2(\vett x) = \nabla \cdot \vett \sigma(\vett x) - \vett f(\vett x)\quad \quad \textrm{in} \quad \quad \Omega
\end{equation}
in which $\vett \sigma(\vett x)$ denotes the Cauchy stress and $\vett f(\vett x)$ is the vector of external forces.

Since a linear elastic and isotropic material is considered here, one can write the residual of the Saint Venant-Kirchhoff law as
\begin{equation}
\vett r_3(\vett x) = \vett \sigma(\vett x) - \mathbb C {:} \vett{\varepsilon}(\vett x) \quad \quad \textrm{in} \quad \quad \Omega
\end{equation}
with $ \mathbb C =  \mathbb C(E,\nu)$ being the spatially constant 4-th order elasticity tensor depending on the material constants $E$ (elasticity modulus) and $\nu$ (Poisson's ratio).
At the boundary the residual of the Cauchy principle is 
\begin{equation}
\vett r_4(\vett x) = \vett t(\vett x) - \vett \sigma(\vett x) \cdot \vett n(\vett x) \quad \quad \textrm{on} \quad \quad \Gamma
\end{equation}
with $\vett u(\vett x)$ as the outward unit normal vector at the boundary of the domain. 
The physical boundary conditions of the underlying Boundary Value problem are stated in their residual form. This leads to 
\begin{equation}
\vett r_5(\vett x) = \vett u(\vett x) - \vett u_0(x) \quad \quad \textrm{on} \quad \quad \Gamma_d
\end{equation}
on the boundary $\Gamma_d$ where the displacement state is given (Dirichlet Boundary Conditions).
The residual of the Neumann boundary condition 
\begin{equation}
\vett r_6(\vett x) = \vett t(\vett x) - \vett t_0(x) \quad \quad \textrm{on} \quad \quad \Gamma_n
\end{equation}
is established on the boundary $\Gamma_n$ with prescribed boundary stress vector state $\vett t_0(\vett x)$.

In the following we decide to fulfil the residuals $\vett r_1$, $\vett r_3$,and $\vett r_4$ exactly, such that
\begin{align}
\vett r_1(\vett x) &= \vett 0 \rightarrow \vett{\varepsilon}(\vett x) = \frac{1}{2} \left(\nabla \vett u(\vett x) + \nabla \vett u^T(\vett x)\right) \quad \quad \textrm{in} \quad \quad \Omega \\
\vett r_3(\vett x) &= \vett 0 \rightarrow \vett \sigma(\vett x) = \mathbb C {:} \vett{\varepsilon}(\vett x)\quad\quad \quad\quad  \quad\quad \quad \quad \,\,\,\, \textrm{in} \quad \quad \Omega \\
\vett r_4(\vett x) &= \vett 0 \rightarrow \vett t(\vett x) = \vett \sigma(\vett x) \cdot \vett n(\vett x)\quad \quad \quad\quad \quad \quad \,\,\,\,\, \textrm{on} \quad \quad \Gamma
\end{align}
and resulting relations can be used directly. On the other hand we assume that the residuals $\vett r_2$ (equilibrium), $\vett r_5$ (Dirichlet boundary conditions) and $\vett r_6$ (Neumann boundary conditions) are not a priori zero and only fulfilled in an approximate sense.

It is possible to train the ANN by defining appropriate cost function which enforces above mentioned governing equations considering the BCs on $\Gamma_d$ and $\Gamma_n$ \citep{SIRIGNANO20181339,Raissi2018}. Let us assume that we have $n_d$ points on $\Gamma_d$, $n_n$ points on $\Gamma_n$, and $n$ points in $\Omega$. The cost function representing the loss of the ANN output is considered as the point-wise evaluated weighted residuals 
\begin{align}
\mathit{MSE}_e = & \frac{1}{n}\sum_{i}^n \vett w_2(\vett x_i) \cdot \vett r_2(\vett x_i) \,\,\, \quad \quad in \quad  \quad \Omega,\label{balance}\\
\mathit{MSE}_d = & \frac{1}{n_d}\sum_{i = 1}^{n_d} \vett w_5(\vett x_i) \cdot \vett r_5(\vett x_i) \quad  \quad \,on \quad  \quad \Gamma_d, \label{EqDBC}\\
\mathit{MSE}_n = & \frac{1}{n_n}\sum_{i = 1}^{n_n} \vett w_6(\vett x_i) \cdot \vett r_6(\vett x_i) \quad  \quad on \quad  \quad \Gamma_n,\label{neumann}\\
\mathit{MSE} = & \mathit{MSE}_d+ \mathit{MSE}_n + \mathit{MSE}_e,\label{eqMSE}
\end{align}
in which the weighting functions are chosen as the respective residual functions
\begin{align}
\vett w_2(\vett x) &= \vett r_2 (\vett x)\\
\vett w_5(\vett x) &= \vett r_5 (\vett x)\\
\vett w_6(\vett x) &= \vett r_6 (\vett x)
\end{align}
and thus leading to a cost function similar to the least-square method, however, they are not based on a space discretisation using a mesh. In other words, one does not need to discretise the whole spatial domain in order to find the solution at a specific point. Basically, the only requirement is to place sufficient data points on $\Gamma_d$ and $\Gamma_n$ which will be studied in details in Section \ref{sec_sensitivity}.
\newpage

Having calculated the cost function, we can choose an optimisation algorithm to minimise it and train the ANN. Here, we develop a method called "hybrid MGA-MSGD" to efficiently obtain a reliable trained ANN which provides us with the continuous solution (displacements and their gradients) of the problem. The abbreviation MGA indicates Modified Genetic Algorithm and MSGD stands for Multi-level Stochastic Gradient Descent.

\section{Hybrid MGA-MSGD training algorithm} \label{sec_MGAMSGD}

 
A Genetic Algorithm (GA) is a stochastic population-based optimisation approach which is considered as a type of Evolutionary Algorithm (EA) \citep{Goldberg1989, Holland1992}. The advantage of this method is that it performs a global search with the expectation to find a basin of attraction of the global minimum of the cost function. However, it is observed that it  is not effective as a local search engine which is the reason for combining it with SGD-based algorithms such as BP \citep{Belew90evolvingnetworks, Yao1999}.  
 
 In this section, we introduce a Modified GA approach enhanced by the Multi-level SGD suitable for both global and local search for better identification of learnable parameters of ANNs. 
 First, we introduce each element of GA and compare it with the ones of MGA in the following.
 
\paragraph{Population} 
  Both algorithms start with constructing a \emph{population}. In GA, a population consists of individuals that each represents a potential solution to a problem while in MGA, the population is a set of individuals that together provide one potential solution. In the presented framework, the population is a list of the binary representation of the network's learnable parameters (weights and biases) so that, for fixed network architecture, they characterise the output of the network. The reason for such arrangement is the high dimension of the learnable parameters and time-consuming procedure of distance function calculation. It also allows us to modify the most important parameters that cause error explosion so that we save computational effort. Furthermore, the binary representation is used so that the strength of mutation, which will be defined, can be controlled and principles of GA can be exploited.

\paragraph{Chromosomes and genes} 
The basic definitions of these elements are the same in GA and MGA. Each individual of a population (learnable parameters of ANN) is called a \emph{chromosome} which consists of several \emph{genes}. The genes can adopt letters, numbers, vectors etc. Here, a chromosome is the binary representation of one learnable parameter, and a gene can adopt one integer of this representation (which is 0 or 1). In other words, a population is a high dimensional array of the binary representation of learnable parameters such as weights and biases. We also consider the sign of the parameter as a gene so that it can be flipped with a given chance maintaining the diversity in the learnable parameters' space. 

\begin{figure}
  \centering
  \includegraphics[width = 12cm]{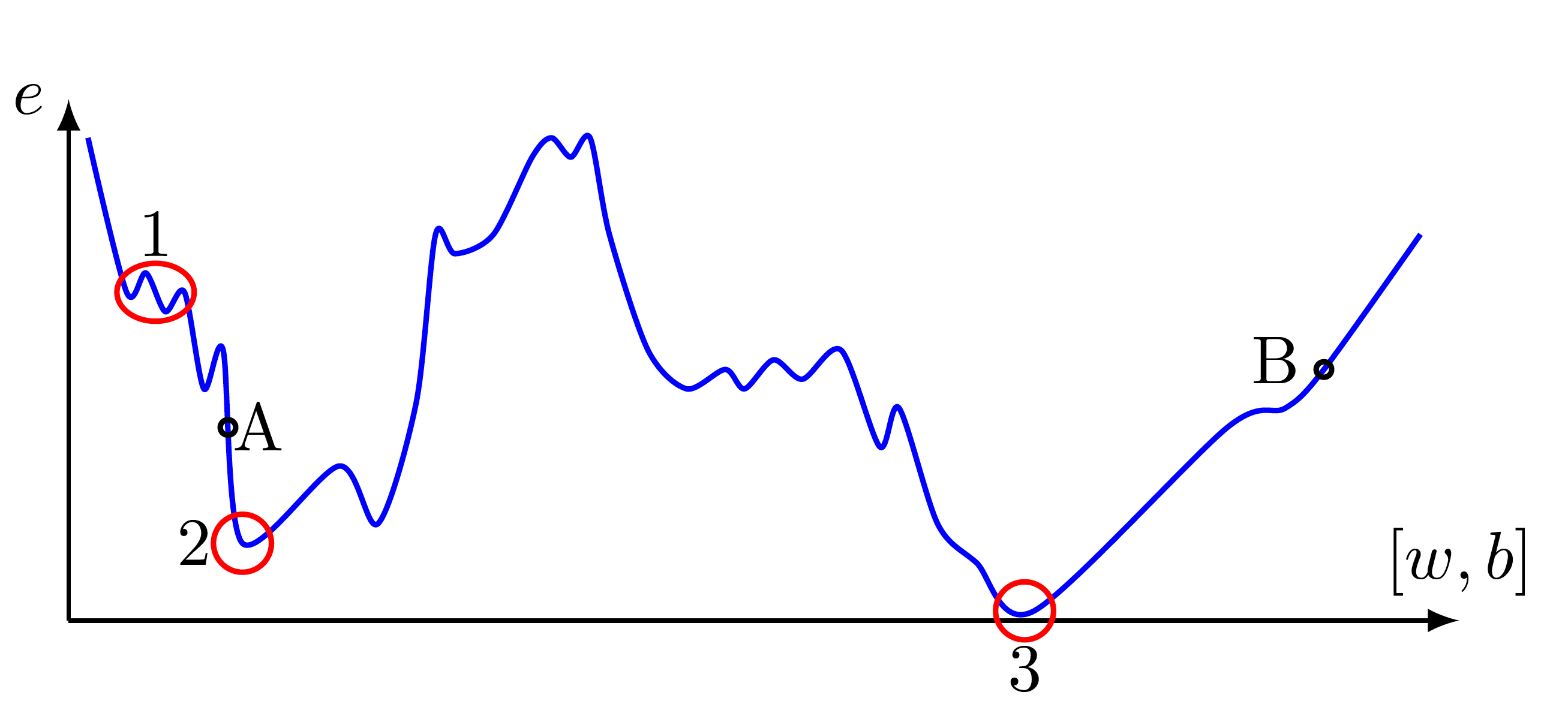}
  \caption{A schematic sketch of the error ($e$) associated with every location in learnable parameter's space $[w,b]$(e.g. weights and biases). Note that the dimension of $[w,b]$ is way higher than what is shown here. Getting trapped in local minima shown in Zone 1 is avoidable choosing a relatively large learning rate that allows larger steps in the horizontal direction. Using SGD, the local minimum in Zone 2 will be the final response even with a large learning rate if we start from its basin of attraction (e.g. Point A). The latter can be tackled by means of MGA which is able to carry out a global path-independent search and step from Point A to Point B. We highlight that although Point A as an initial location has less error, the Point B is favourable as it is in the neighbourhood of "global" minimum.}
  \label{fig_LocalM}
  \end{figure}

\begin{figure}
  \centering
  \includegraphics[width = 8cm]{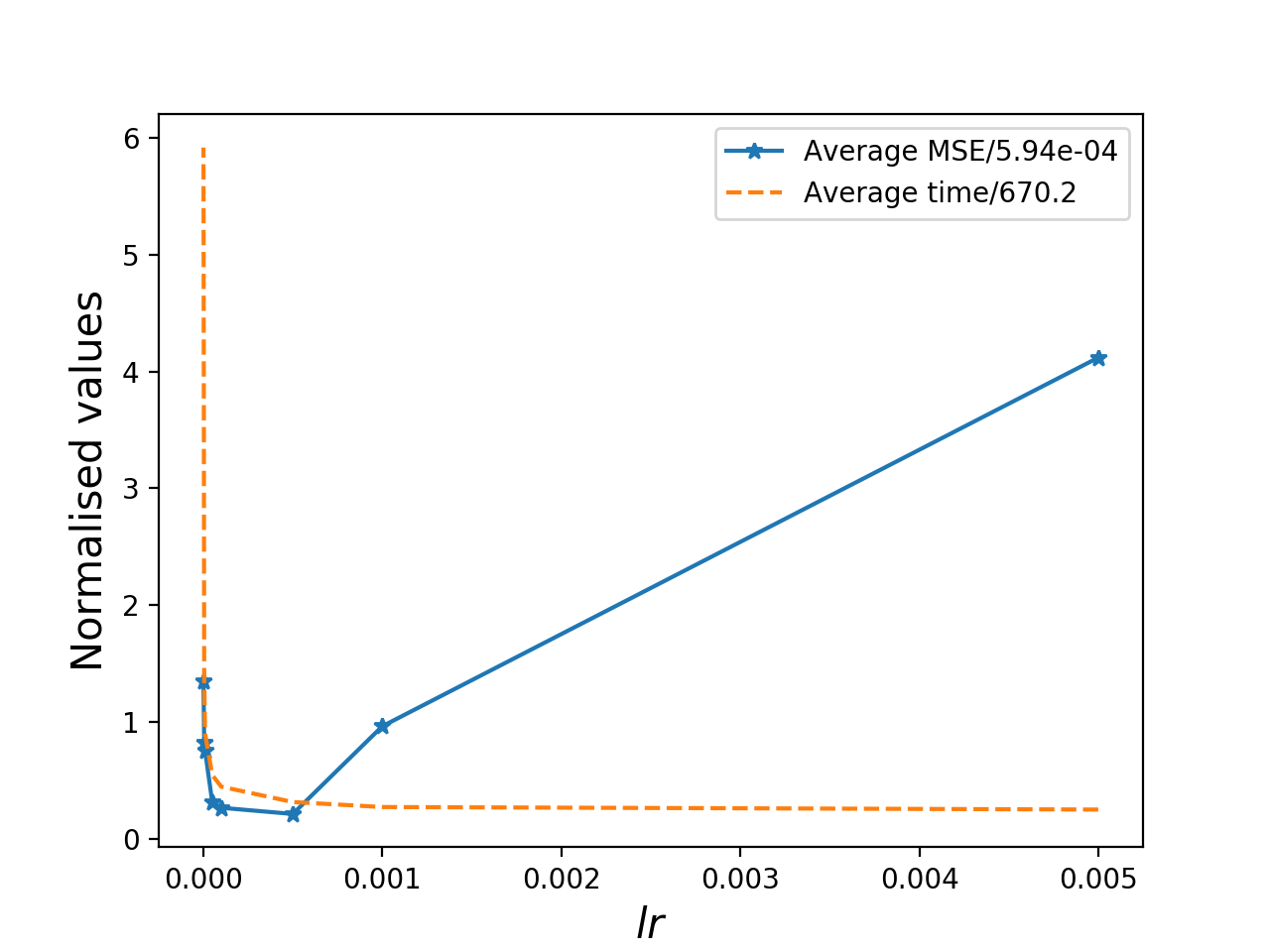}
  \caption{The profile of learning rates vs average $\mathit{MSE}$ shows that having smaller learning rates does not necessarily result in  higher accuracy. The average $\mathit{MSE}$ (or other variables in the plots) is its averaged value of several analyses under equivalent conditions so that decreasing the effect of randomness and indicating the role of the parameter under study.
   This plot is produced experimentally by performing training tests with different learning rates. Every point is the average of the value of interest ($\mathit{MSE}$ and time) of 100 tests.}
  \label{fig_lrF}
  \end{figure}

\paragraph{Fitness and Importance} 
In GA, \emph{fitness} is a value assigned to each chromosome which is the relative success of the corresponding individual on the problem. For example, if we choose this value between 0 and 100, the individual that provides the correct value has the fitness 100, and the one that has no similarity with the correct one adopts 0.  Due to the fact that in MGA the individuals do not provide a solution to the problem (the entire population can only provide it) they can not be directly examined and assigned a fitness value. However, using gradient-based decomposition approaches such as Back Propagation we can quantify the relative sensitivity of the cost function with respect to the learnable parameters (i.e. $\frac{\partial C}{\partial \sbt}$, where $\sbt$ is one learnable parameter and $C$ is the cost function) showing how "important" is each individual. We highlight that the most important individuals are the ones that cause error explosion enforcing the choice of small learning rates. The latter can be seen as one reason of getting trapped into local minima which are avoidable using larger learning rates, e.g. the ones in Zone 1 of Figure \ref{fig_LocalM}. This non-trivial dependency of the results on the learning rate is proven via sensitivity analysis and shown in Figure \ref{fig_lrF}.
 Consequently, we replace the term "fitness" in GA with the term "importance" in MGA.

\paragraph{Selection and offsprings} 
In GA, some chromosomes (according to their fitness scores) are chosen in the process of \emph{selection}. Usually, the chromosomes with higher fitness score have a higher chance of being selected. These chromosomes are called "parents" and are chosen to produce offsprings from which the next generation is constructed which is to be more similar to the selected chromosomes while maintaining the diversity.  
In MGA, the chromosomes to be modified are selected based on their importance via Tournament Selection procedure (\cite{Miller1995GeneticAT}). The selected parameters for minimisation of the cost function (the most important chromosomes), as explained before, are the ones that cause problems such as error explosion enforcing the choice of very small learning rate which, in turn, increases the chance of getting trapped into a local minimum. The fact that GA tends to modify the parameters with high sensitivity (importance) helps us to avoid this problem. These individuals, subsequently, go through some operations out of which we hope to acquire a better arrangement of the population. We highlight that, here, the unselected individuals are the surviving ones.

 Using the activation function $ELU$ in Equation \eqref{activeF} and expanding the BP gradient formula $\frac{\partial C}{\partial \sbt}$ we understand that the sensitivity of a parameter does not only depend on itself. In fact, it depends on all the parameters that are directly or indirectly in contact with that parameter. Consequently, in order to maintain the diversity of selected individuals (i.e. avoiding the similar choices in several MGA iterations), we introduce the condition that MGA can not select the same parameters more than a specific times (here, twice). This restriction is released once all of the parameters in the population are previously selected (in order to avoid empty selection) by clearing the history of the selected individuals.

\paragraph{Operations} 
 There are two well-known operations of GA called \emph{crossover} and \emph{mutation}. The idea of the former is to produce some offsprings by combining the genes of two parents in order to generate new individuals that provide new solutions. The produced offsprings are partially similar to their parents. However, when we use MGA for ANN training, there is no need for the similarity between the selected individuals (learnable parameters) as the inheritance from the previous population is obtained via the unselected individuals tuned only by the connectionist optimiser gradient-based CSGD. So, as expected and from our observations as well as the literature \citep{Yao1999, Heimes1997, Sarkar1997, BORNHOLDT1992327}, the crossover does not perform well in ANN training. On the other hand, the second operation (mutation) is very useful in the global/local search for the optimal point as it maintains the genetic diversity by flipping the value of one or more genes \citep{Goldberg1989}. Here, we also introduce the sign of the chromosome as a gene that can be flipped by a given chance. The sign and the first following genes are responsible for the global search while the last genes are responsible for small changes performing a local search.  In other words, if we flip the first/last genes of the chromosome, a large/small change is imposed thus a global/local modification on the selected parameters takes place. We introduce three factors controlling global/large scale, medium scale, or local/small scale mutation probability. The first one (global mutation) allows us to, for example, step from Point A to Point B in Figure \ref{fig_LocalM} while the others are responsible for smaller steps. 

\begin{figure}
  \centering
  \includegraphics[width = 8cm]{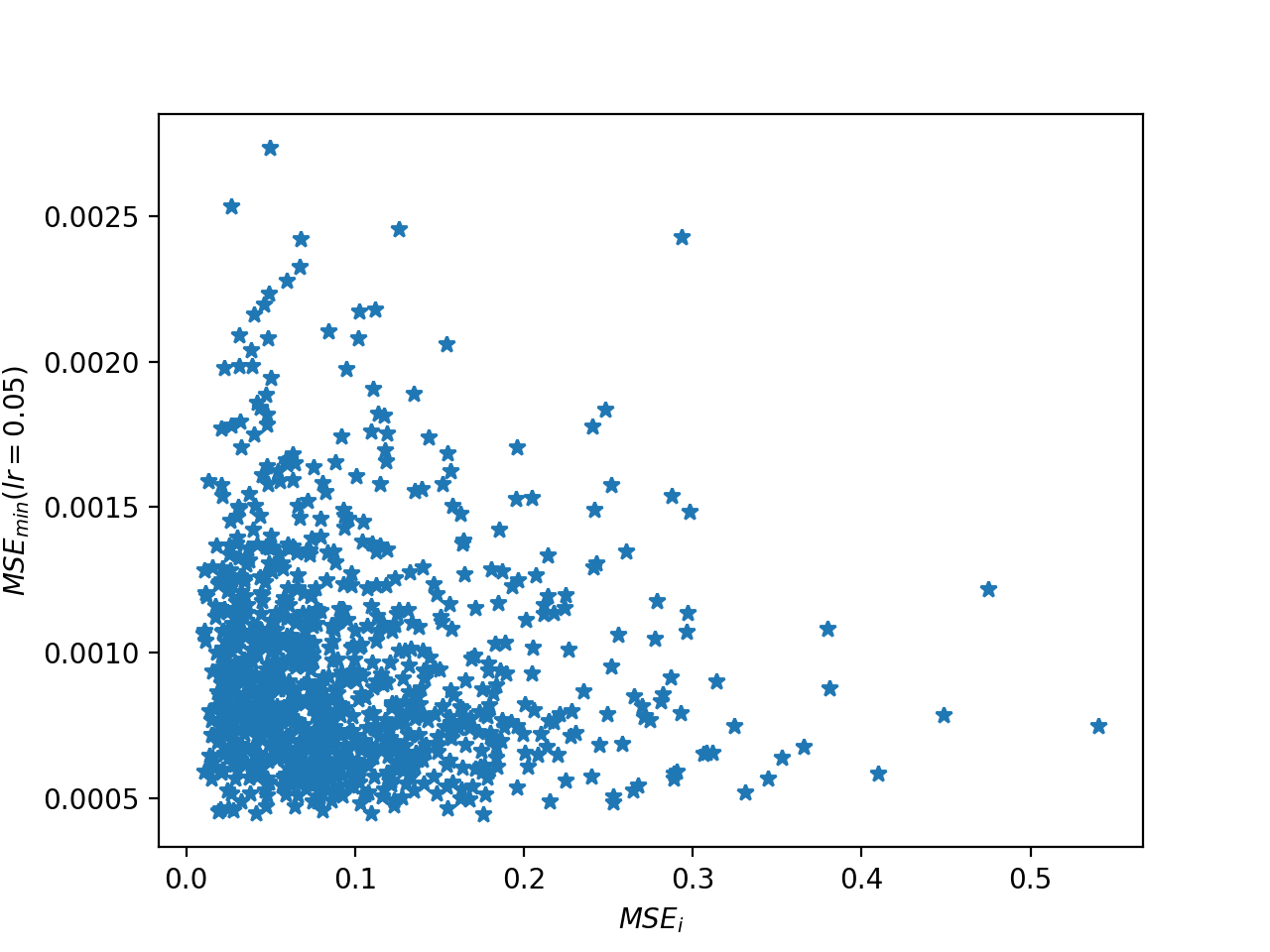}
  \caption{The plot of $\mathit{MSE}_{min}$ vs $\mathit{MSE}_i$ shows that although $\mathit{MSE}_i$ can be an important parameter, it is not sufficient as the only measure to identify a "better" population.}
  \label{msecomp}
  \end{figure}

 
 \paragraph{Offspring, population qualification, and new generation} 
  The obtained chromosomes from mutation are, then, converted back to the float representation constructing the \emph{offspring}. The latter is the link to produce a new population which we hope to fall into a better point in the parameters' space. In fact, we replace the selected parameters by their corresponding offspring with the expectation that this new arrangement is "better" than the previous one. In this case, a "better" population is the one that leads us nearer to the global minimum (accuracy) at fewer increments and time (efficiency). The new population serves as the initial parameters for a SGD optimisation process. At this point, the important question is that if a population that creates a smaller initial cost function $\mathit{MSE}_i$ is a better one or it is possible that we achieve, after performing SGD, a more accurate final result (shown by $\mathit{MSE}_{min}$) from a larger $\mathit{MSE}_i$? Figure \ref{msecomp} shows that having smaller $\mathit{MSE}_i$ does not necessarily result in a smaller $\mathit{MSE}_{min}$, which is obtained via performing a SGD optimisation with a fine learning rate to find the nearby minimum. 

\begin{figure}
  \centering
  \begin{subfigure}{6.6cm}
  \includegraphics[width = 6.7cm]{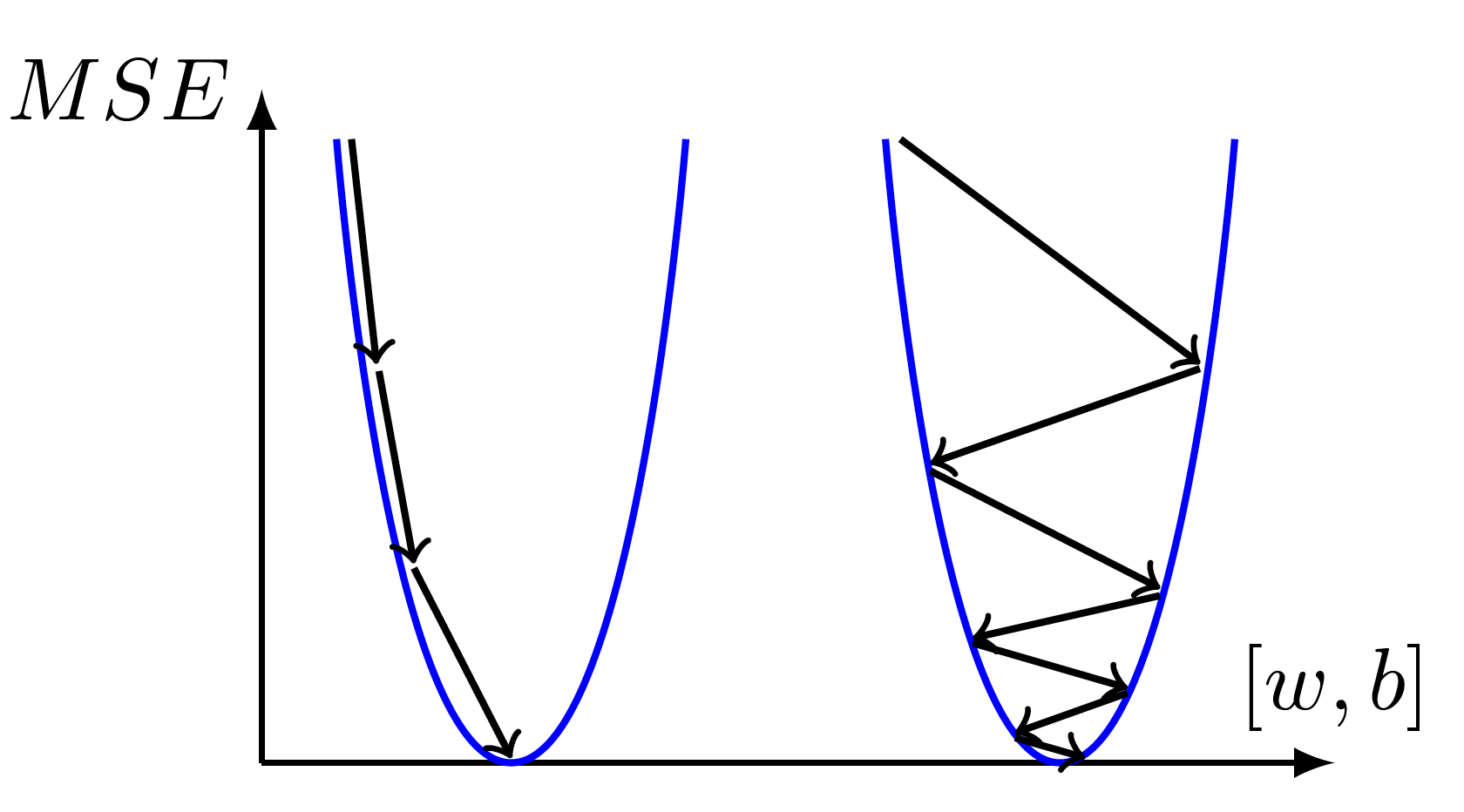}
  \caption{Two conditions under which the cost function converges to a minimum with sufficiently small learning rates.}
  \label{fig_conv}
  \end{subfigure}\quad
  \begin{subfigure}{4.8cm}
  \includegraphics[width = 5.1cm]{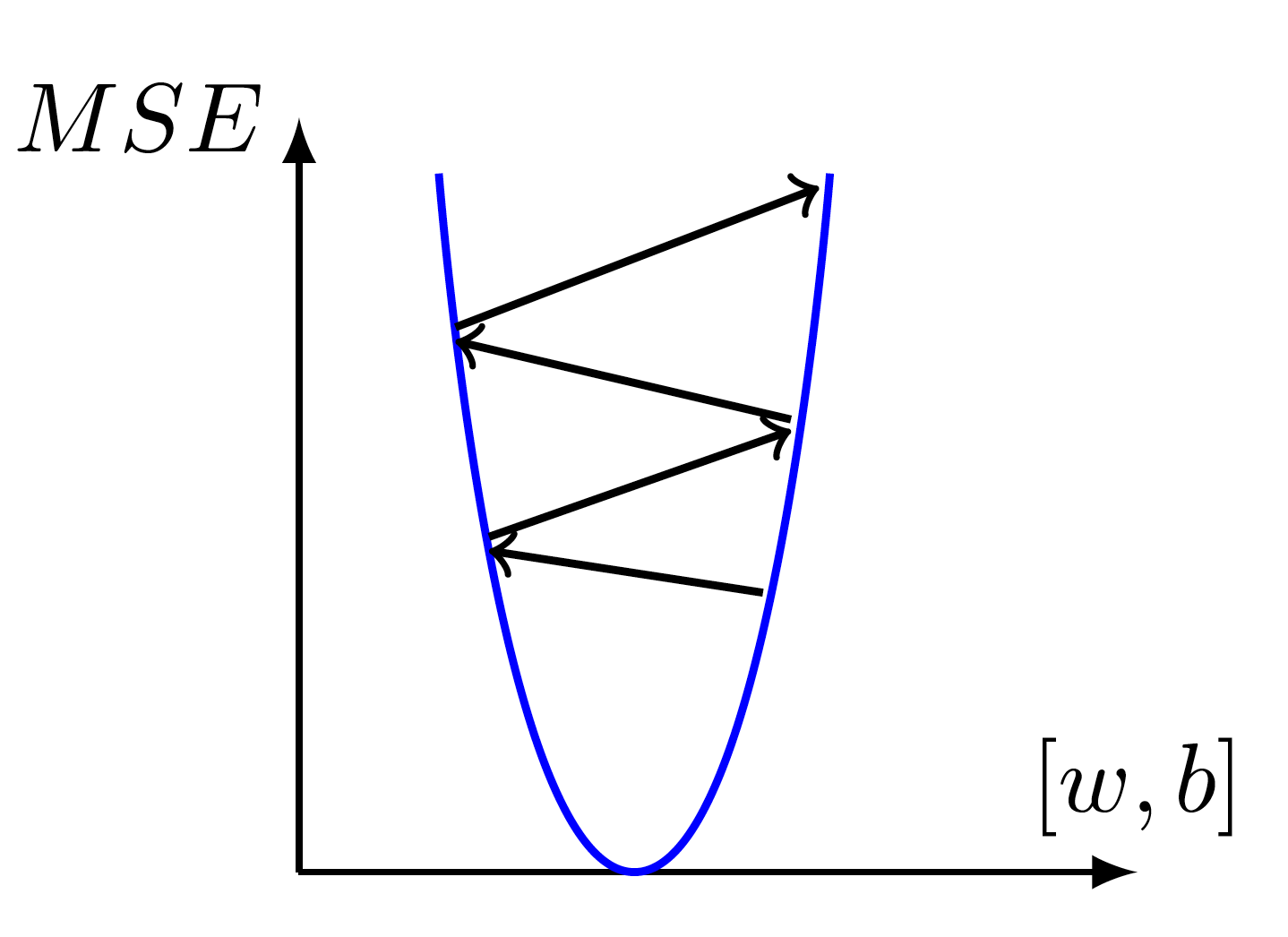}
  \caption{The divergence of the cost function from the minimum due to a large learning rate (large horizontal step).}
  \label{fig_div}
  \end{subfigure}
  \caption{A schematic representation of convergence and divergence of the cost function. Although a large learning rate helps to skip several local minima, it causes divergence at a distance from the target minimum. To tackle this condition, we introduce a fine level SGD after achieving the final MGA-CSGD result.}
  \label{fig_div_conv}
  \end{figure}

At this stage, we can introduce the \emph{Coarse-level} SGD (CSGD or a SGD with a large learning rate) so that we can both examine the new population and get near the objective accuracy. In general, if we have a good arrangement of the learnable parameters such that the output/cost function does not largely depend on some specific small group of parameters with very high gradients ($\frac{\partial C}{\partial \sbt}$) it is possible to reach an accurate result with large learning rate and small number of iterations avoiding local minima in a short time. In fact, achieving such population arrangement, we use the full potential of our ANN in which all the neurones play a considerable role in the final outcome. Consequently, we perform a CSGD optimisation based on the new population. Then, we can base the population qualification on $\mathit{MSE}_c$ as follows
\begin{numcases}
{if} \mathit{MSE}_c^{(j)}<\mathit{MSE}_c^{(j-1)}:  \quad\textrm{qualified population} \\
\mathit{MSE}_c^{(j)}>\mathit{MSE}_c^{(j-1)}: \quad \textrm{unqualified population},
\end{numcases}
where $j$ is the current MGA iteration number. If the population gets qualified we call the obtained arrangement (the one after CSGD) a new \emph{generation} and base the following MGA iterations on it otherwise we neglect it and we use the previous population/generation. Due to the complexity of the BP (using the third order partial derivatives) we might need to impose a break condition when the cost/distance divergence starts in order to avoid stepping away from the minimum  \citep{Duckett2008}. The latter condition, as well as the convergence one, is shown in Figure \ref{fig_div_conv}.

Finally, at the end of the MGA analysis with CSGD based qualification, we perform a \emph{Fine-level} SGD (FSGD or a SGD with small learning rate) in order to achieve a more accurate result with smaller $\mathit{MSE}$. The reason is shown in Figure \ref{fig_div_conv}. This completes the hybrid MGA-MSGD training procedure. 

\section{Sensitivity analysis and numerical results} \label{sec_sensitivity}
The above-described framework has several non-mechanical/framework parameters that may considerably affect the efficiency of the procedure and accuracy of the results. Given the parametric uncertainty, it seems necessary to employ a Sensitivity Analysis method (SA) to identify important model parameters and to calibrate them. Furthermore, due to the stochastic nature of the ANN training procedure, every piece of data provided for the sensitivity analysis (e.g. mean $\mathit{MSE}$, mean time, standard deviation, variance, etc.) should be obtained by averaging the results of several analyses which brings high computational expense and should be considered as a constraint when adopting a SA method. Therefore, the Morris method \citep{Morris1991}, which is a simple and efficient approach based on changing one-factor-at-a-time (OAT) is employed.

\begin{figure}
	\includegraphics[width=7cm]{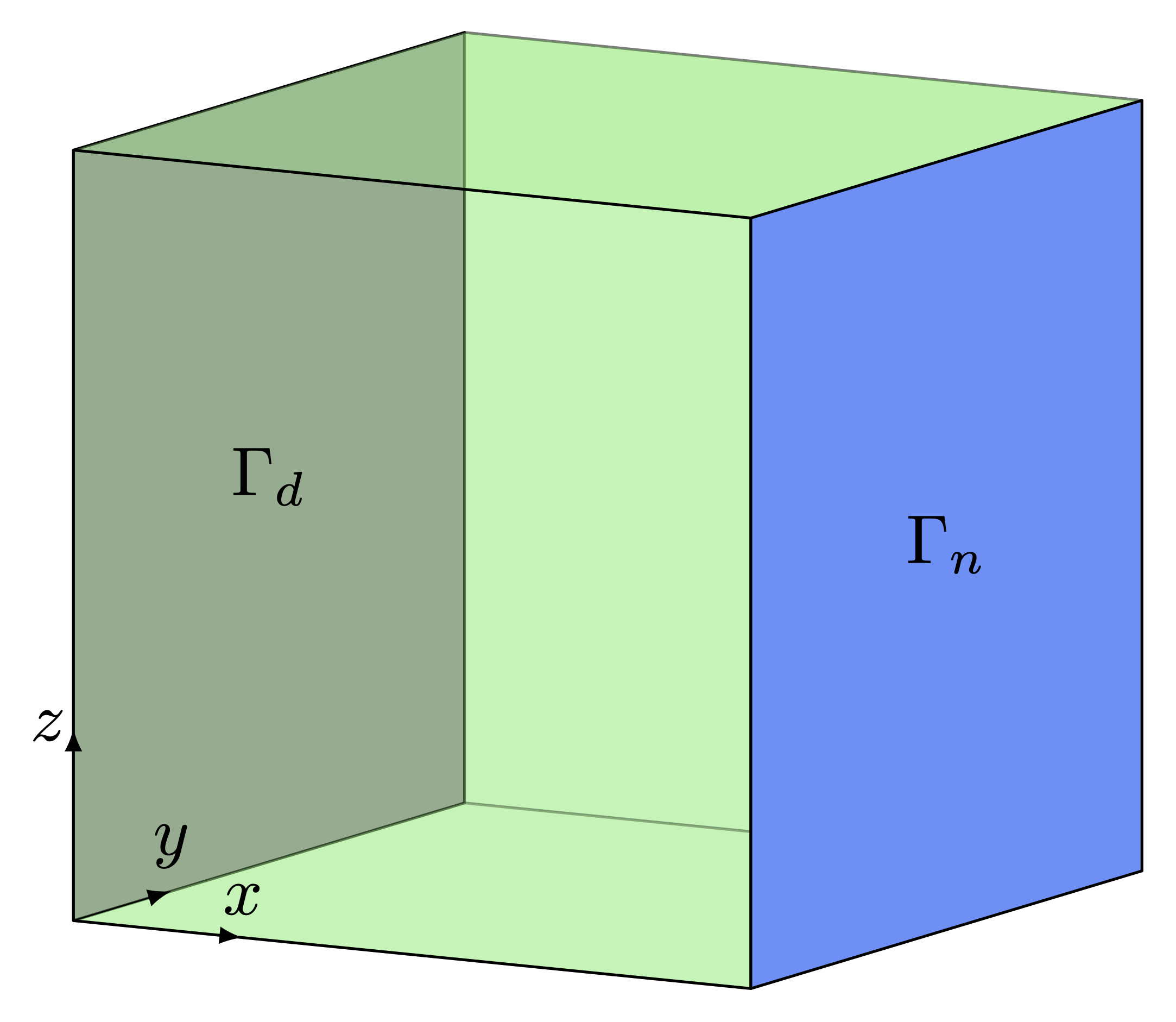}
	\centering
	\caption{A representative sketch of the model. Zero displacements in $\vett x$ direction ($u_{x0}=0$) as the Dirichlet BC is applied on $\Gamma_d$ and $\vett t = (-0.1, 0, 0)$ Neumann BC is applied on $\Gamma_n$.}
	\label{fig_cube}
\end{figure}

\paragraph{\textbf{Statement of the test problem}}
The deformation behaviour a unit cube of homogeneous isotropic linear elastic material with zero displacements in $\vett x$ direction ($u_{x0} = 0$) on $\Gamma_d$ located at $x=0$ and a nonzero uniformly distributed $\vett t_0 = (-0.1, 0, 0)$ (see  Figure \ref{fig_cube}).  The uniqueness conditions
\begin{align}
\frac{1}{n} \sum_{i=1}^n  u_y = 0 \label{uniquex} \\
\frac{1}{n} \sum_{i=1}^n  u_z = 0 \label{uniquey}
\end{align}
are also imposed. 

The density and distribution of the training dataset can be controlled by the number of the sampling points in $\vett x$, $\vett y$, and $\vett z$ directions, which are shown by $N_{\vett x}$, $N_{\vett y}$, and $N_{\vett z}$, and the number of points located on the boundaries with Dirichlet/Neumann constraint $N_b$  to the one of the whole volume (the total number) $N_V$. We also introduce the variable $\beta_i$ that modifies the value of $\frac{N_b}{N_V}$
\begin{equation}
\frac{N_b}{N_V} = \frac{N_{bu} + \beta_i N_{bu}}{N_{Vu} + \beta_i N_{bu} }
\end{equation}
Where $N_{Vu}$ and $N_{bu}$ are the total number of the points and the number of the ones exactly on the boundaries at the uniform distribution of sampling points, repectively. Note that, as this method is free of space discretisation (mesh-free), one can randomly choose the sampling points up to the condition that there must be some points on the surfaces $\Gamma_d$ and $\Gamma_n$ so that we can enforce Dirichlet/Neumann BCs in the training procedure. However, we choose to have an entirely controlled point distribution in order to understand how it affects the outcome of the analysis. 

 In Morris method there are two parameters showing the sensitivity of the output to a variable which are the standard deviation $\sigma$ and mean elementary effect $\mu$ computed by
 \begin{align}
 \sigma = \sqrt{\frac{\sum_{i=1}^{n}(x_i - \bar x)^2}{n-1}}\\
 \mu =  \frac{ \sum_{i=1}^{n}\lvert(x_i - \bar x) \rvert}{n-1} 
 \end{align}
 Generally, the more the value of $\sigma$ and $\mu$ the more important the corresponding parameter.
 
 We study the effects of the non-mechanical/framework parameters shown and explained in Table \ref{MorrisSensitivity}, where surviving population fraction is the fraction of unselected population which are the ones that remain unchanged in the next population. The results are provided in Table \ref{MorrisSensitivity1}.

\begin{table}
\caption{}
\centering
\scalebox{0.85}{
\begin{tabular}{|c| c| c| c|}
\hline
Parameter & Notation & Distribution & Base value\\
\hline
Coarse scale learning rate & $lr_c$ & $U(0.5, 1)$ & $0.7$\\
\hline
Number of MGA iterations & $N_{GAi}$ & $U(10,60)$ & $30$\\
\hline
Number of ANN hidden layers & $N_h$ & $U(2,6)$ & $3$ \\
\hline
Number of neurones in each hidden layer & $N_{nh}$ &  $U(3, 20)$ & $10$\\
\hline
Surviving population fraction & $P_{sf}$ &  $U(0.9, \approx 1)$ & $0.98$\\
\hline
Number of sampling points in $\vett x$ direction & $N_{x}$ &  $U(5, 14)$ & $10$ \\
\hline
$ N_{\vett x} - N_{\vett y} =  N_{\vett x} - N_{\vett z}$ & $N_{x} - N_y$ & $U(0, 8)$ &   $0$ \\
\hline
Sampling point distribution parameter & $\beta_i$ & $U(0, 2)$ &  $0$\\
\hline
Global mutation probability & $M_g$ & $U(0.1, 0.5)$ & $0.1$ \\
\hline
Median mutation probability & $M_m$ & $U(0.1, 0.5)$ & $0.1$ \\
\hline
Local mutation probability & $M_l$ & $U(0.1, 0.5)$ & $0.1$ \\
\hline
\end{tabular}}
\label{MorrisSensitivity}
\end{table}

\begin{table}
\caption{}
\centering
\scalebox{0.85}{
\begin{tabular}{|c| c| c| c| c| c| c| c| c|}
\hline
\multirow{2}{*}{Notation} & \multicolumn{2}{|c|}{Average $\mathit{MSE}$} & \multicolumn{2}{|c|}{Minimum $\mathit{MSE}$} & \multicolumn{2}{|c|}{Average time} & \multicolumn{2}{|c|}{Minimum time} \\ \cline{2-9}
& $\mu$ & $\sigma$ & $\mu$ & $\sigma$  & $\mu$ & $\sigma$  & $\mu$ & $\sigma$ \\
\hline
$lr_c$ &  $2.75e{-06}$ & $9.14e{-06} $  & $8.35e{-07}$ & $2.64e{-06}$  & $9.84$ & $31.1$ & $1.73$ & $5.47$\\
\hline
$N_{GAi}$ & $2.21e{-06}$ & $7e{-06}$ & $1.07e{-06}$ & $3.38e{-06}$ & $13.6$ & $43.16$ & $5.9$ & $18.66$ \\
\hline
 $N_h$ & $2.36e{-05}$ & $4.71e{-05}$ & $9.14e{-06}$ & $1.82e{-05}$ & $10.6$ & $21.34$ & $4.1$ & $8.17$ \\
\hline
$N_{nh}$ & $7.17e{-06}$ & $1.9e{-05} $ & $1.3e{-06}$ & $3.4e{-06} $ & $8.67$ & $22.9$ & $12.58$ & $33.3$ \\
\hline
$P_{sf}$ & $1.92e{-06}$ & $6.4e{-06}$ & $8.8e{-07}$ & $2.86e{-06}$ & $3.36$ & $11.13 $ & $2.22$ & $7.37$\\
\hline
$N_{x}$ & $1.34e{-06}$ & $4.1e{-06}$ & $1.34e{-06}$ & $4e{-06}$ & $30.9$ & $92.8$ & $6.08$ & $18.2$\\
\hline
$N_{x} - N_y$ & $8.18e{-07}$ & $2.31e{-06} $ & $1.24e{-06}$ & $3.5e{-06} $ & $20.29$ & $57.4$ & $4.59$ & $13$\\
\hline
$\beta_i$ & $1.1e{-06}$ & $3.32e{-06}$ & $1.18e{-06}$ & $3.54e{-06}$ & $21.7$ & $65.1$ & $4.19$ & $12.56$\\
\hline
$M_g$ & $8.8e{-07}$ & $1.76e{-06}$ & $7.48e{-07}$ & $1.5e{-06}$ & $11.5$ & $23$ & $19.6$ & $39.3$\\
\hline
$M_m$ & $1.51e{-06}$ & $3.03e{-06}$ & $7.07e{-07}$ & $1.41e{-06}$ & $59.4$ & $118.8$ & $52.3$ & $104.7$\\
\hline
$M_l$ & $5.53e{-07}$ & $1.1e{-06}$ & $5.64e{-07}$ & $1.13e{-06}$ & $39.5$ & $79$ & $20.1$ & $40.3$\\
\hline
\end{tabular}}
\label{MorrisSensitivity1}
\end{table}

\begin{figure}
\centering
	\begin{subfigure}{5.8cm}
	\includegraphics[width=5.8cm]{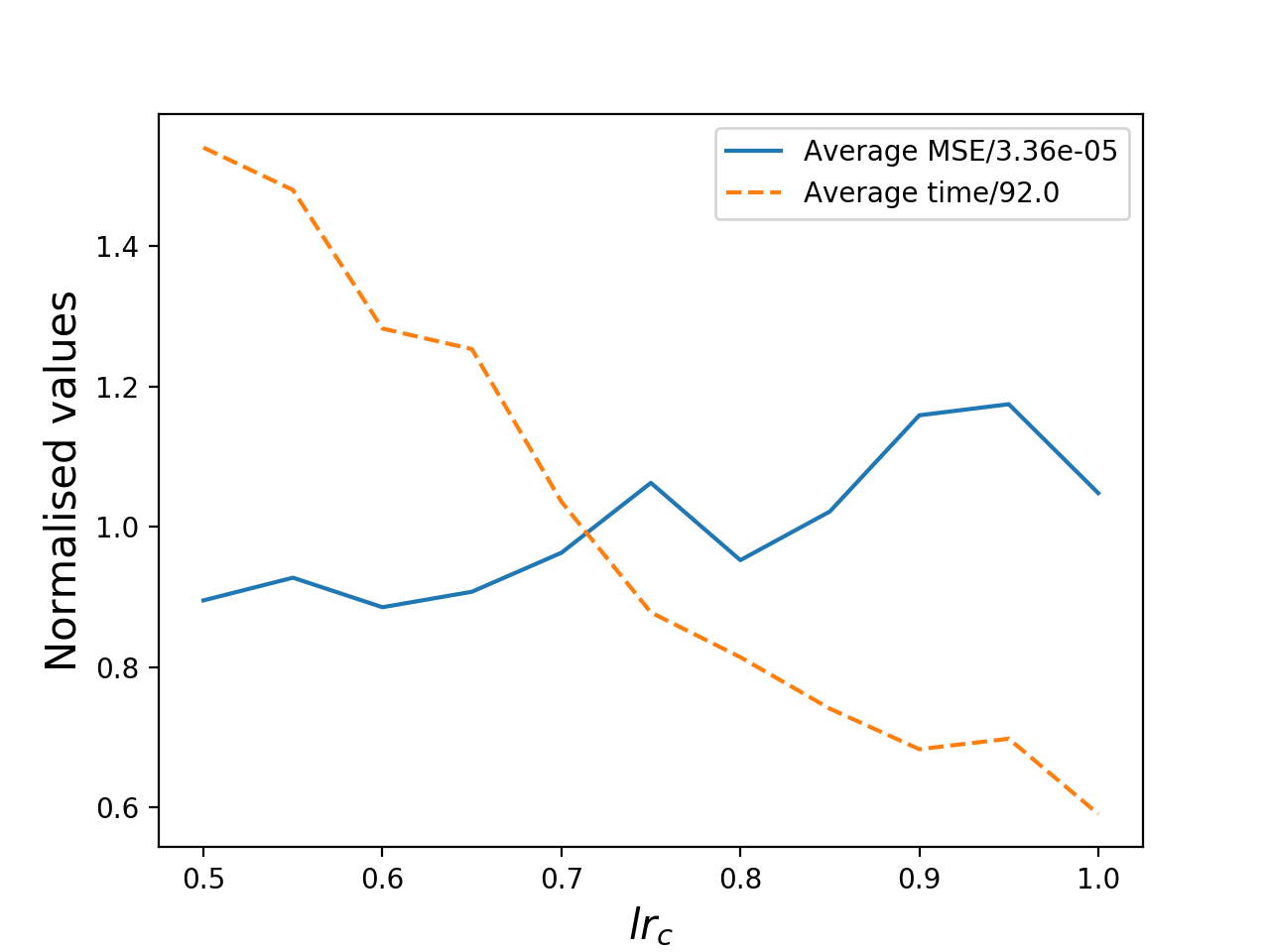}
	\centering
	\caption{Coarse level learning rate ($lr_c$) has a significantly larger impact on the efficiency compared with the accuracy.}
	\label{LRcAVG}
	\end{subfigure} \quad
	\begin{subfigure}{5.8cm}
	\includegraphics[width=5.8cm]{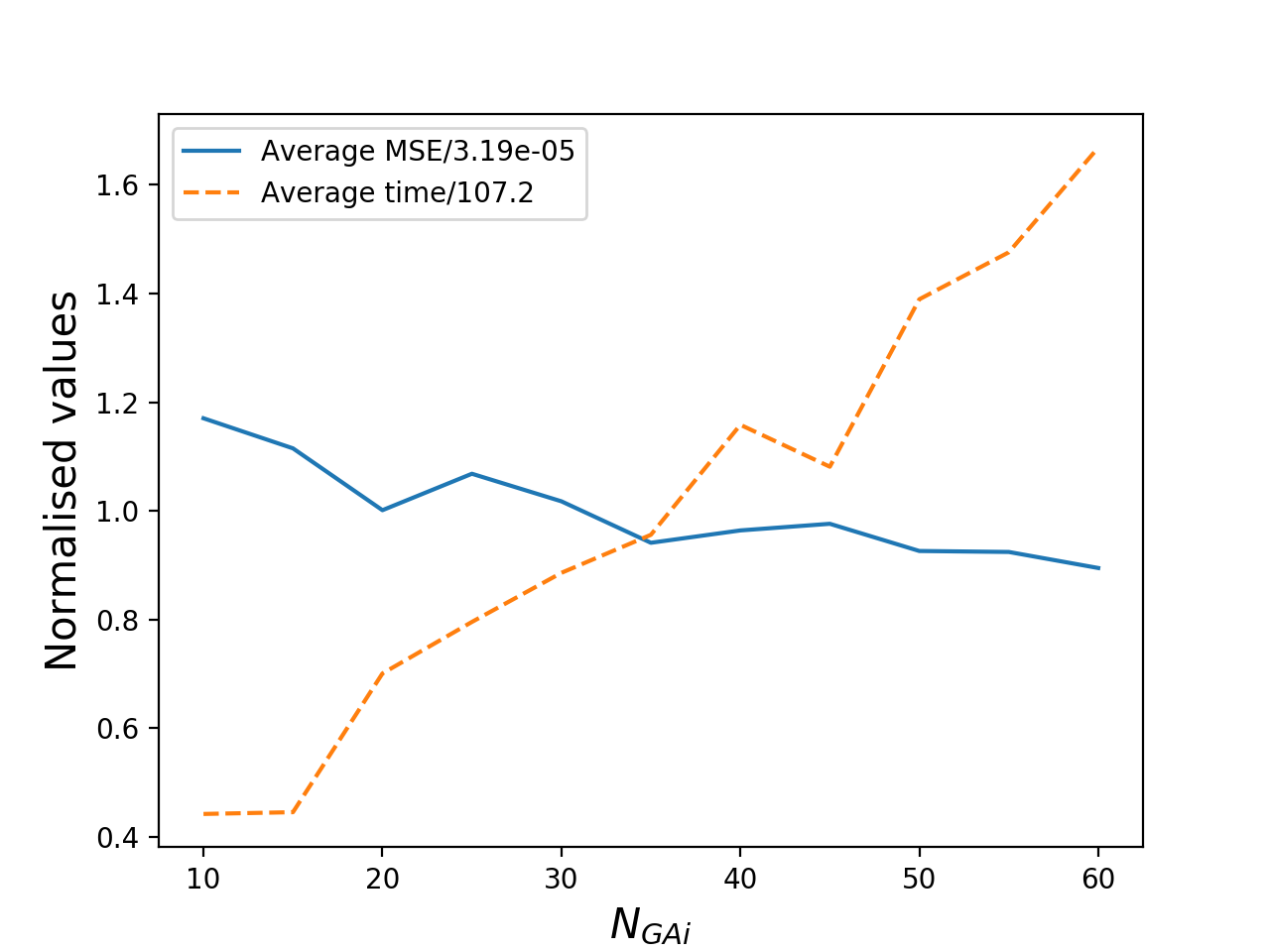}
	\centering
	\caption{The number of iterations of MGA-CSGD ($N_{GAi}$) also has a significantly larger impact on the efficiency compared with the accuracy.}
	\label{NGAiAVG}
	\end{subfigure} 
	\caption{The dependency of the output accuracy and training time on $lr_c$ and $N_{GAi}$ for parameter tuning purpose. The chosen values in this study are $lr_c=0.6$ and $N_{GAi}=30$. Each one of the "average $\mathit{MSE}$" and "average time" values is divided by a specific value shown in the labels which are their mean values (the mean "average time" and "average $\mathit{MSE}$").}	
	\label{propss}
\end{figure}

\begin{figure}
\centering
	\begin{subfigure}{5.8cm}
	\includegraphics[width=5.8cm]{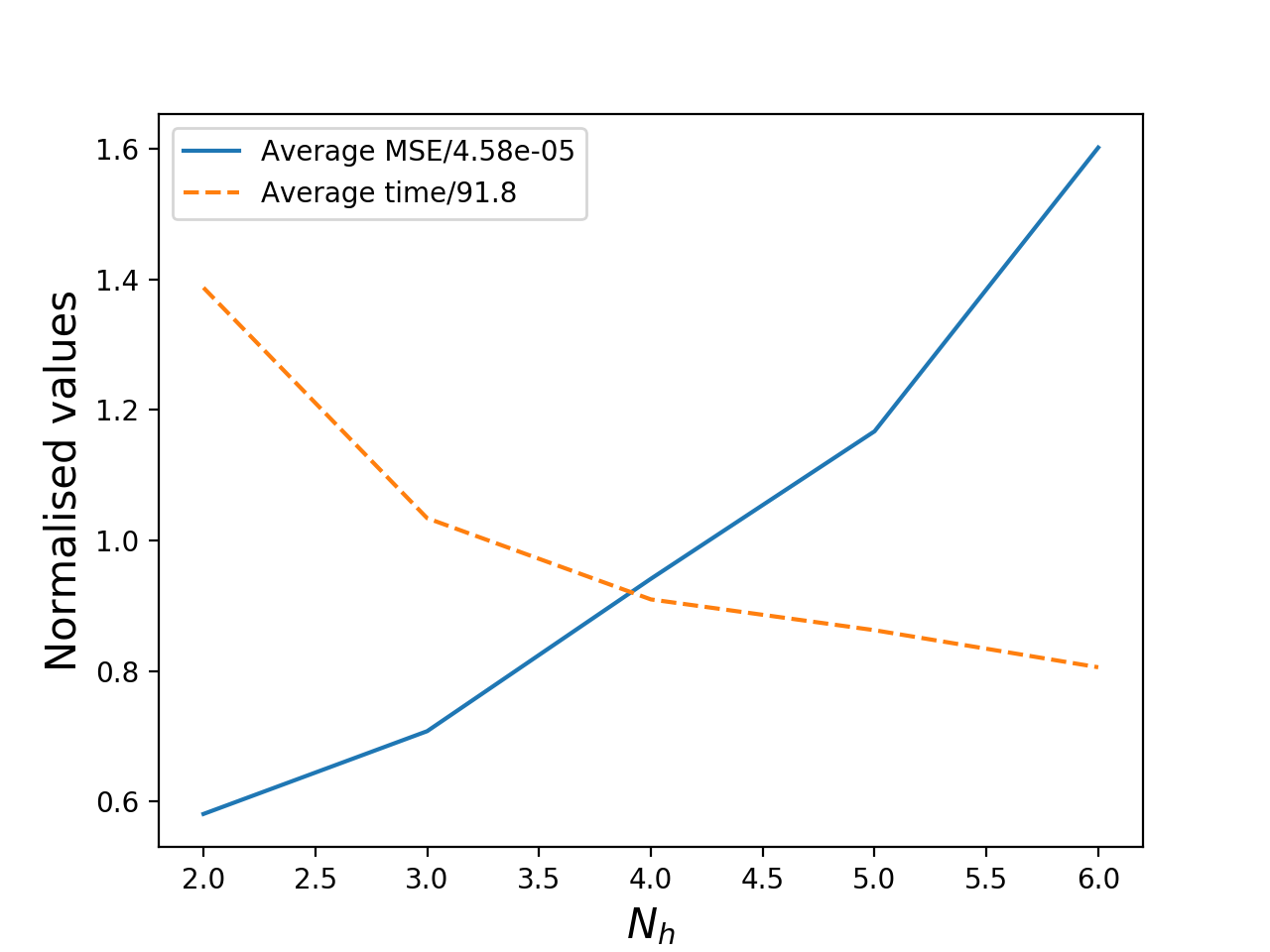}
	\centering
	\caption{The highest dependency of the result's accuracy is on the number of hidden layers indicating that a network with two hidden layers produces the lowest $\mathit{MSE}$ at an affordable time which agrees with the observations in \citep{Dockhorn2019}. Note that deeper networks can also produce better results, but they need more MGA-CSGD iterations etc.}
	\vspace{0.6cm}
	\label{NhAVG}
	\end{subfigure} \quad 
	\begin{subfigure}{5.8cm}
	\includegraphics[width=5.8cm]{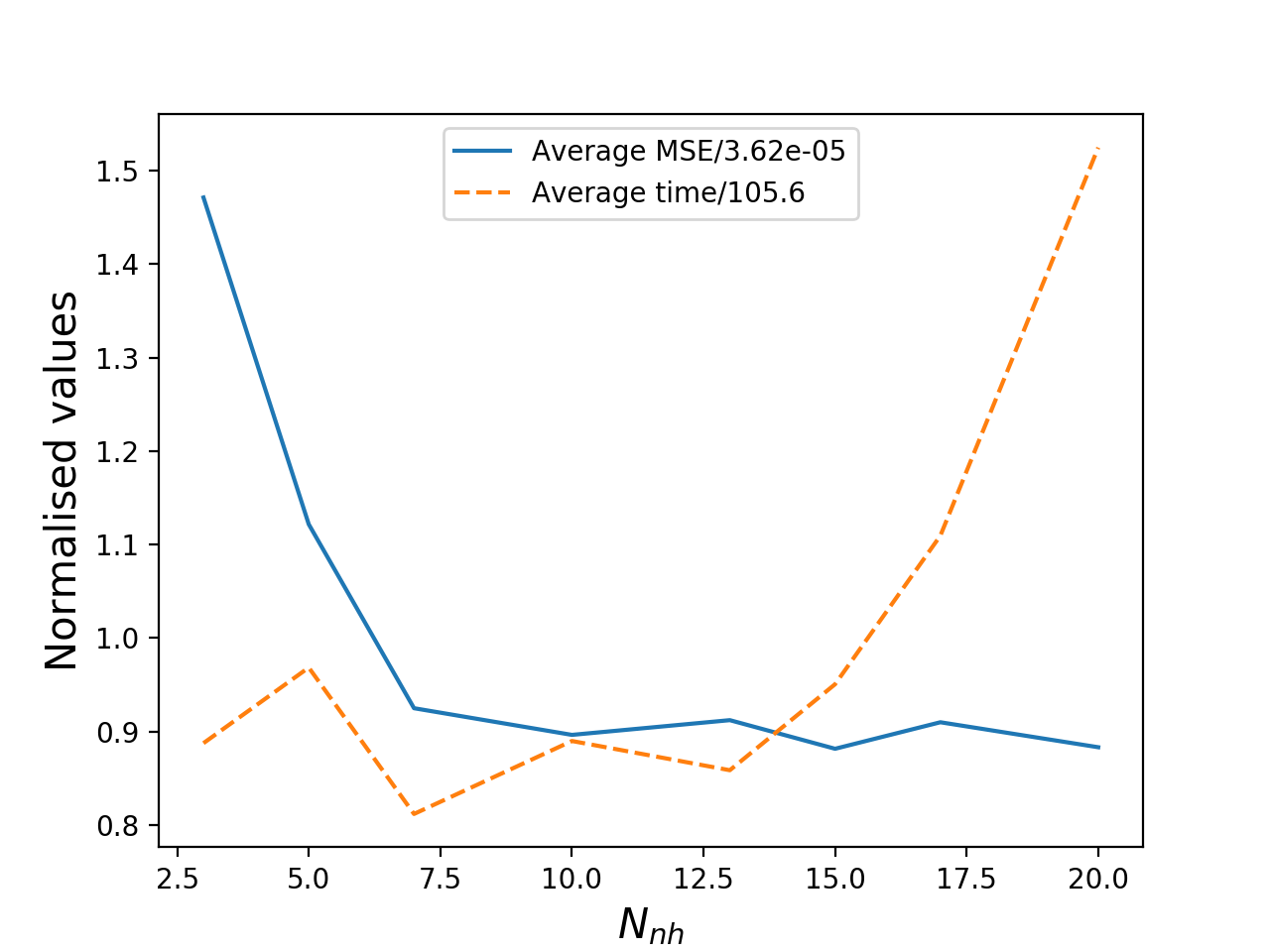}
	\centering
	\caption{The profile of the number of neurones in each hidden layer ($N_{nh}$) vs normalised $\mathit{MSE}$ and time indicates that before a specific point the accuracy of the results depends on the number of the neurones to a great extent which decreases sharply after it. In fact, this point can be interpreted as the minimum required number of neurones that are able to interpolate the specific problem. For this problem, we choose $N_{nh}=10$.}
	\label{NnhAVG}
	\end{subfigure} 
	\caption{The number of hidden layers and neurones in each are of the most critical parameters to be tuned.  Obtaining an acceptable accuracy at such a small number of neurones and layers demonstrates the efficiency of the presented method.}
\end{figure}

\begin{figure}
\centering
	\begin{subfigure}{5.8cm}
	\includegraphics[width=5.8cm]{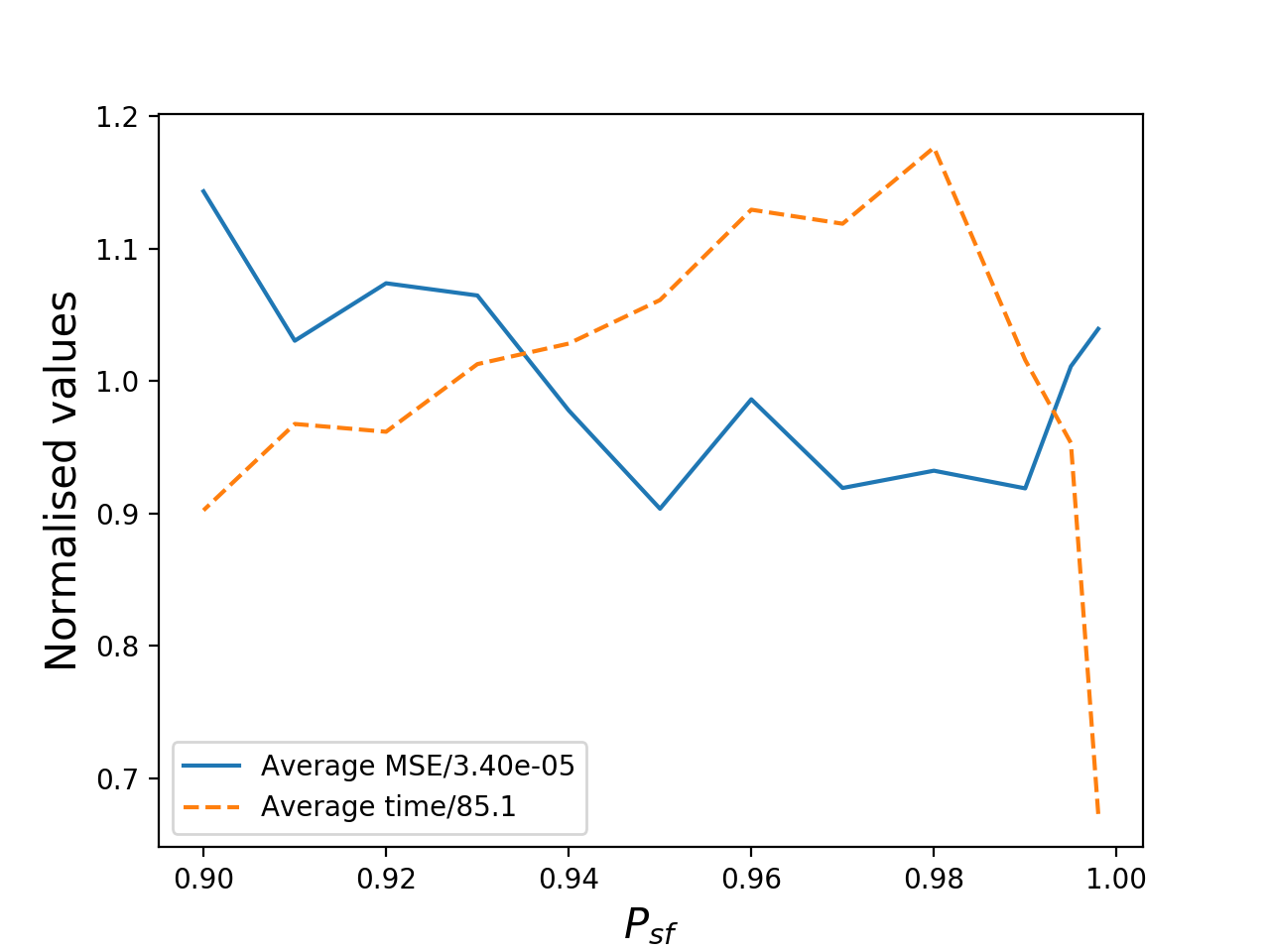}
	\centering
	\caption{Surviving population fraction $P_{sf}$ is considerably larger than the ones in typical applications of GA as the target here are only those demanding small learning rates. As it is shown here, average $\mathit{MSE}$ decreases and average time increases at increasing $P_{sf}$ up to 0.99.}
	\vspace{1.13cm}
	\label{PsfAVG}
	\end{subfigure} \quad 
	\begin{subfigure}{5.8cm}
	\includegraphics[width=5.8cm]{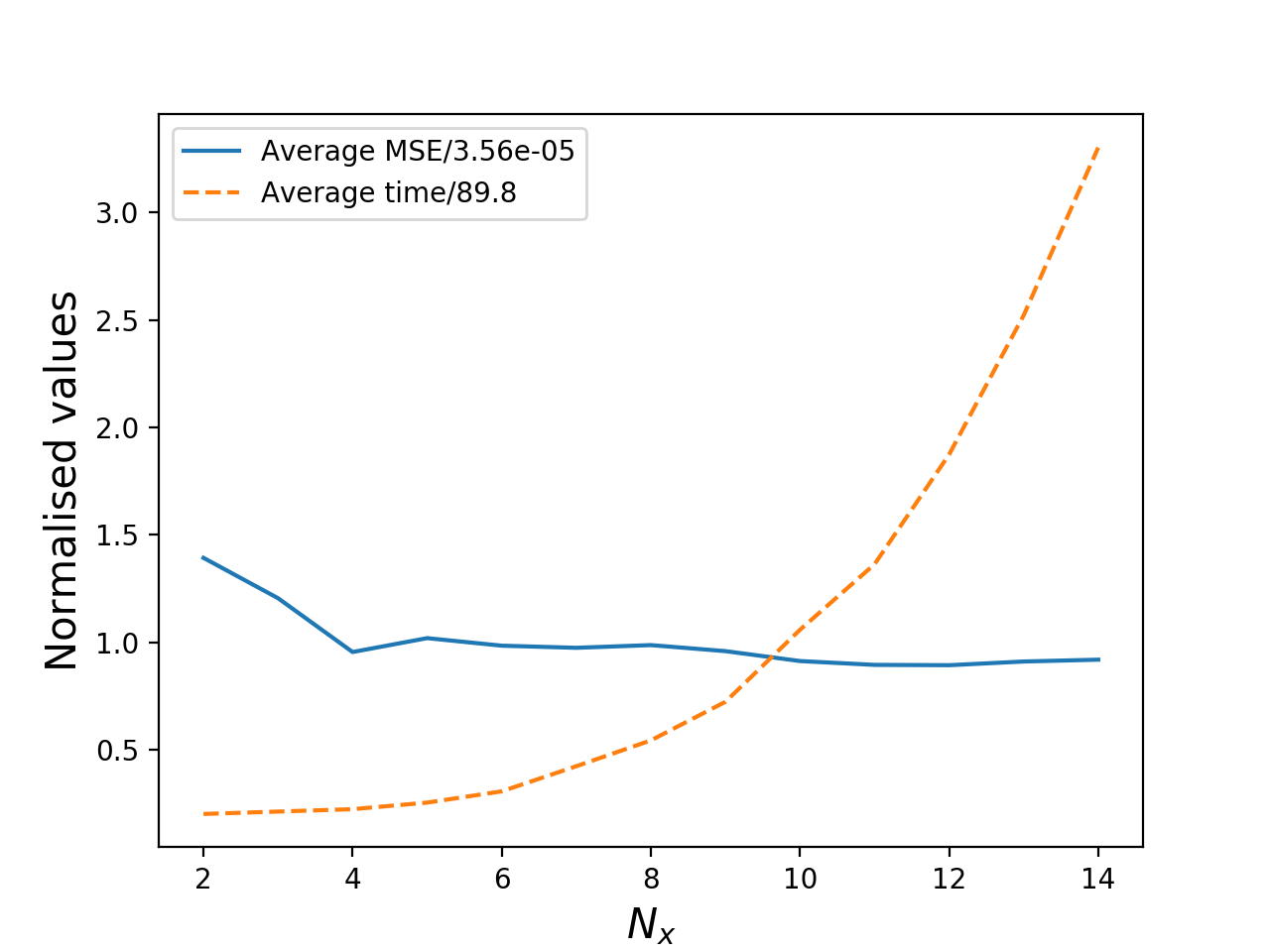}
	\centering
	\caption{$N_x$ represents the density of the training/sampling points. It has a profile similar to the one of the number of neurones in each hidden layer in terms of average $\mathit{MSE}$ meaning that it should be more than a specific number from where on the accuracy is not highly sensitive on it. This can be seen as one significant achievement of the proposed approach as increasing it further the required time increases exponentially due to the gradient calculation.}
	\label{NxAVG}
	\end{subfigure} 
	\caption{The profile of the surviving population fraction and the density of training points vs training time and obtained $\mathit{MSE}$ as a measure of accuracy.  Each one of the "average $\mathit{MSE}$" and "average time" values is divided by a specific value shown in the labels which are their mean values (the mean "average time" and "average $\mathit{MSE}$").}
\end{figure}

\begin{figure}
\centering
	\begin{subfigure}{5.8cm}
	\includegraphics[width=5.8cm]{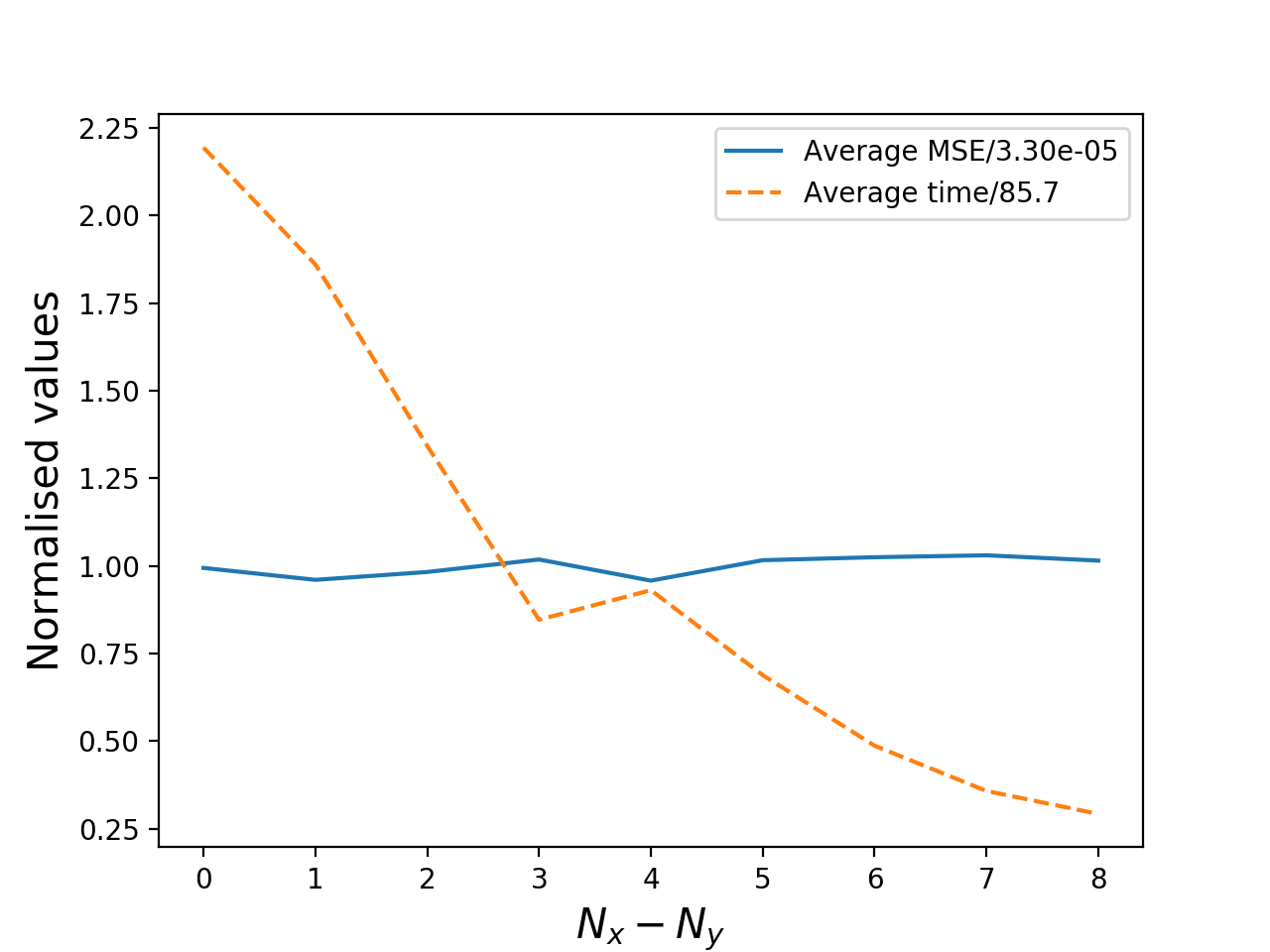}
	\centering
	\caption{The profile of $N_x-N_y$ represents the training/sampling points' distribution in different directions showing a relatively low dependency of the obtained result on it. Note that $N_x=10$ is fixed at its base value so that the decreasing time at increasing $N_x-N_y$ is due to the lower number of training points.}
	\label{Nx-NyAVG}
	\end{subfigure} \quad 
	\begin{subfigure}{5.8cm}
	\includegraphics[width=5.8cm]{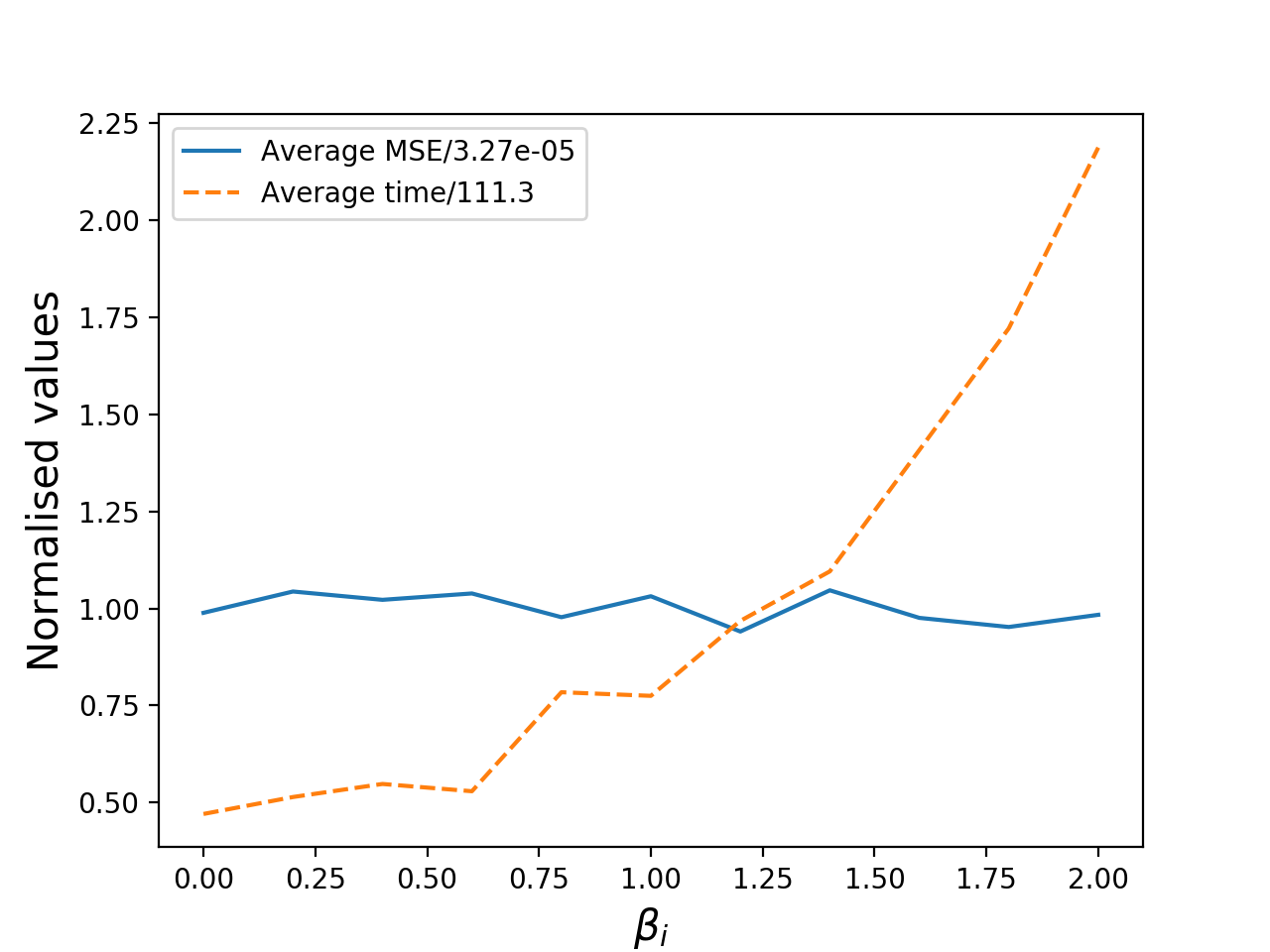}
	\centering
	\caption{$\beta_i$ determines the distribution of training points in a way that increasing it the number of the points exactly on the boundaries increase.}
	\vspace{1.15cm}
	\label{betaAVG}
	\end{subfigure}
	\caption{The small sensitivity of the results on the two factors ($N_x-N_y$ and $\beta_i$) controlling the distribution of the training points results in the choice of uniform distribution.}
	\label{Nx-Ny_beta}
\end{figure}

\begin{figure}
\centering
	\begin{subfigure}{5.8cm}
	\includegraphics[width=5.8cm]{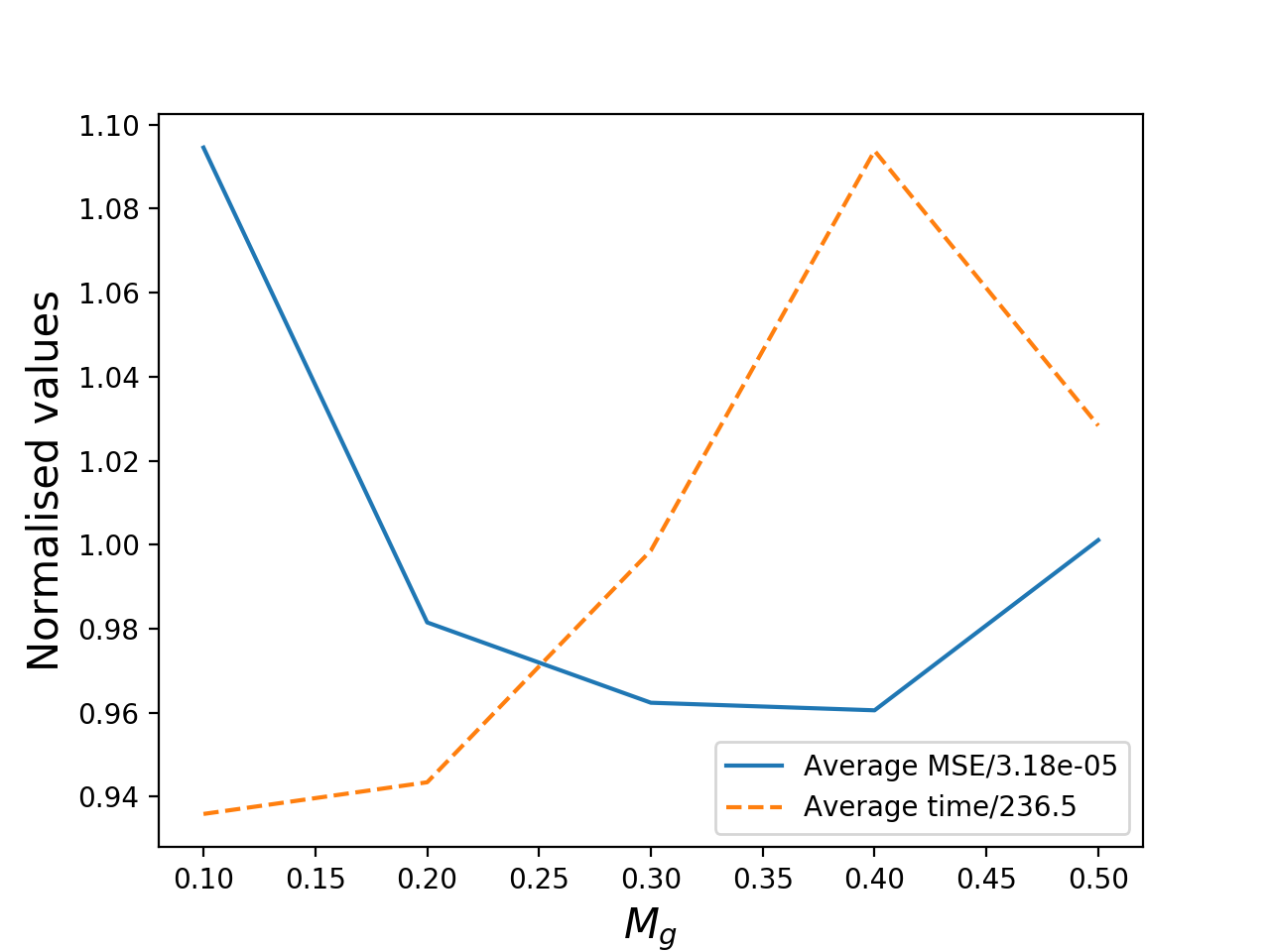}
	\centering
	\caption{Global mutation probability around $M_g=0.3$ provides acceptable results. This factor is responsible for large scale steps in the learnable parameters' space.}
	\label{MgAVG}
	\end{subfigure} \quad 
	\begin{subfigure}{5.8cm}
	\includegraphics[width=5.8cm]{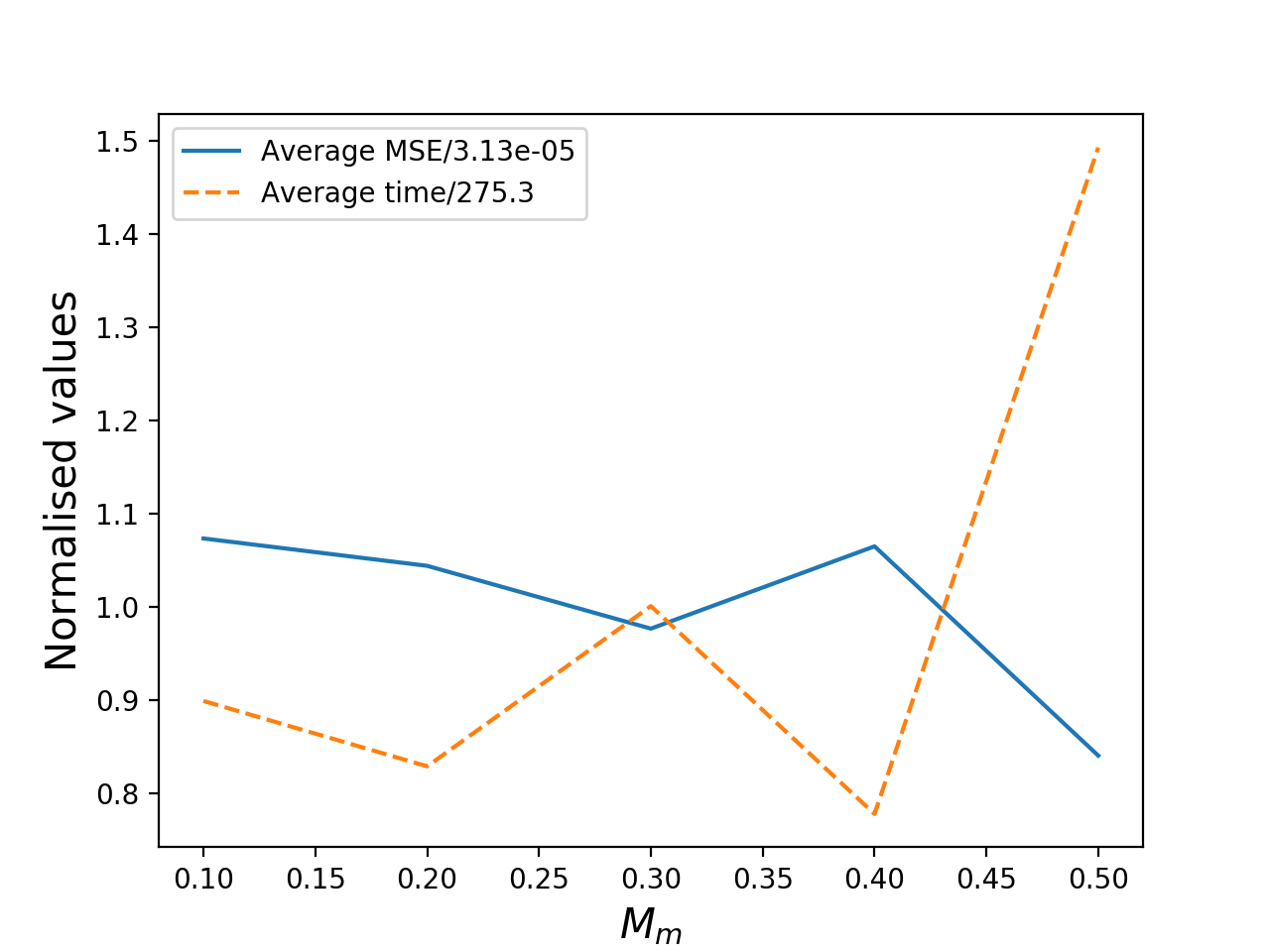}
	\centering
	\caption{$M_m$ controls the median steps in ANN training. In this case, a probability of 0.3 is chosen.}
	\vspace{0.6cm}
	\label{MmAVG}
	\end{subfigure}
	\label{MgMm}
	\begin{subfigure}{5.8cm}
	\includegraphics[width=5.8cm]{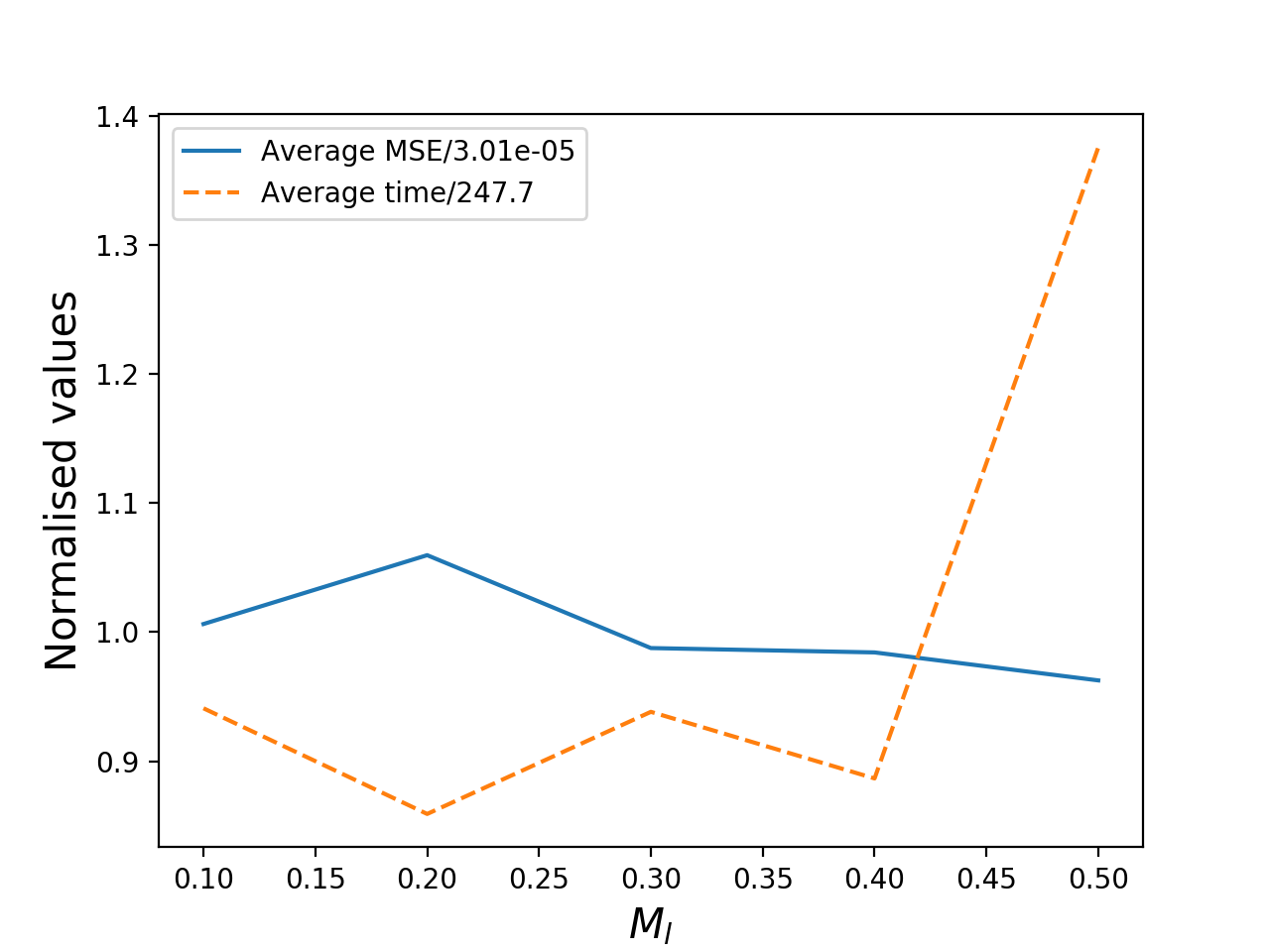}
	\centering
	\caption{A small path-independent step together with SGD (path-dependent) guarantees that we take the full advantage of the local search. $M_l=0.3$ is shown to be an acceptable value.}
	\label{MlAVG}
	\end{subfigure}
	\caption{The probability of mutation for large, medium, and small steps has a considerable effect on the training result.}
	\label{props}
\end{figure}

We plot the normalised values of average $\mathit{MSE}$ and training time with respect to their mean value (which is provided in the line labels) in Figures \ref{propss}-\ref{props} in order to tune the algorithm's parameters. In general, the training time depends on several factors including the dimension of population, number of generations, number of sampling points etc. Although the plots in the provided figures are somehow noisy, one can understand the overall trends. The goal is to decrease $\mathit{MSE}$ at the lowest possible cost/increase in time. For example, according to Figure \ref{NGAiAVG}, although by decreasing the coarse level learning rate ($lr_c$) the average $\mathit{MSE}$ decreases, the training time increases sharply, so we choose a moderate value such as $lr_c =0.6$. The inverse of the latter profile, somehow, holds for the number of MGA-CSGD iterations shown in Figure \ref{NGAiAVG}. We again choose a median value for this parameter $N_{GAi} = 30$. The most important parameters that also have the greatest values of Morris factors for average $\mathit{MSE}$ are the number of hidden layers $N_h$ and neurones in each layer $N_{nh}$. By increasing the former, as shown in Figure \ref{NhAVG}, the accuracy of the results decreases showing that the provided training method performs the best with a small number of hidden layers. It is noteworthy that the decrease in average time at increasing $N_h$ is due to the fact that we have fewer generations (qualified populations) when we have more hidden layers. 
Usually, the loss divergence in unqualified populations happens at a low number of CSGD iterations (sometimes at the first try) indicating the advantage that we mostly spend time on suitable arrangements of the learnable parameters which are most probable to result in high accuracy. 

Figure \ref{NnhAVG} shows that, at least, seven neurones are required in each hidden layer for an acceptable accuracy ($N_{nh}=10$ is chosen). In fact, this shows that using hybrid MGA-MSGD training method we exploit a high potential of the assigned hidden layers and neurones. Figure \ref{PsfAVG} shows a noisy profile and indicates that the surviving population fraction $P_{sf}=0.97$ is an appropriate choice. The profile of $N_x$, $N_x - N_y$, and $\beta_i$ that are, respectively, shown in Figures \ref{NxAVG}, \ref{Nx-NyAVG}, and \ref{betaAVG} demonstrate the low sensitivity of average $\mathit{MSE}$ to the number and distribution of sampling points. This is a crucial advantage of the provided method when ANN is applied to solving mechanical problems because, as expected, the average training time grows exponentially at growing number of sampling points. Note that, as in Figure \ref{Nx-NyAVG} the number of $N_x$ is fixed at its base value, by changing the value of $N_x-N_y$ both distribution and number of sampling points are changed.


 For more detailed analysis of sampling point distribution, we consider two different distribution scenarios Case 1: $(N_x, N_y, N_z) = (30,2,2)$ and Case 2: $(N_x, N_y, N_z) = (5,5,5)$ with their results ( average $\mathit{MSE}$, minimum $\mathit{MSE}$, average time, and minimum time)  provided in Table \ref{tab:T2}. The results show that for a smooth mechanical problem, the choice of uniform distribution of data points provides the best $\mathit{MSE}$ and training time. On the other hand, the small difference between the results of Case 1 and Case 2 together with Figure \ref{Nx-Ny_beta}, reveals that being mesh-free and directly differentiable, this approach (solving mechanics problem via ANN) is able to provide an approximate response of specific coordinates independent of additional points inside the body provided that the BCs are appropriately enforced. For example, Case 1 has points only on four sides of the cube, not inside the volume. This feature together with the flexibility and straight forward implementation and providing continuous spatial gradients can be seen as advantages of approximating mechanical problems via ANNs.

\begin{table}
\caption{}
\centering
\scalebox{0.8}{
\begin{tabular}{|c| c| c| c| c|}
\hline
sampling point distribution & Average $\mathit{MSE}$ & Minimum $\mathit{MSE}$ & Average time & Minimum time \\
\hline
$(N_x, N_y, N_z) = (30,2,2)$ & $3.35e{-05}$ & $2.2e{-05}$ & $47$ &  $24.1$\\
\hline
$(N_x, N_y, N_z) = (5,5,5)$ & $3.18e{-05}$ & $2.15e-{05}$ & $38.2$ & $18.9$\\
\hline
\end{tabular}}
\label{tab:T2}
\end{table}

In summary, the chosen setting for this problem is $(N_x, N_y, N_z) = (5,5,5)$, $lr_c=0.6$, $N_{GAi}=30$, $N_h=2$, $N_{nh} = 10$, $P_{sf} = 0.97$, $\beta_i=0$, and $lr_f = 1e{-5}$. To authors' experience, the parameters $lr_c$, $(N_x, N_y, N_z)$, and $N_{nh}$ are the ones that should be available at a higher programming level to the users for tuning purposes.
In the following subsection the results of the presented methodology are compared with the standard training approaches.
\subsection{Training performance}
The most important points to be considered for computational methods are efficiency and accuracy. Here, we choose two standard well celebrated training algorithms, namely classical SGD  (\cite{Robbins2007ASA}) and Adam (\cite{kingma2014adam}), to assess the relative success of the presented framework. Here, for the sake of clarification and reader's convenience, the update strategy for SGD and Adam is provided. The former has simple update terms as below
\begin{equation}
\sbt^{(i)}_{kj}:=\sbt^{(i)}_{kj}- n \frac{\partial \mathit{MSE}}{\partial \sbt^{(i)}_{kj}} 
\end{equation}
where $n$ is a coefficient that determine the learning rate (size of steps) and $\sbt$ is the updated parameter.

\begin{figure}
\centering
	\begin{subfigure}{5.8cm}
	\includegraphics[width=5.8cm]{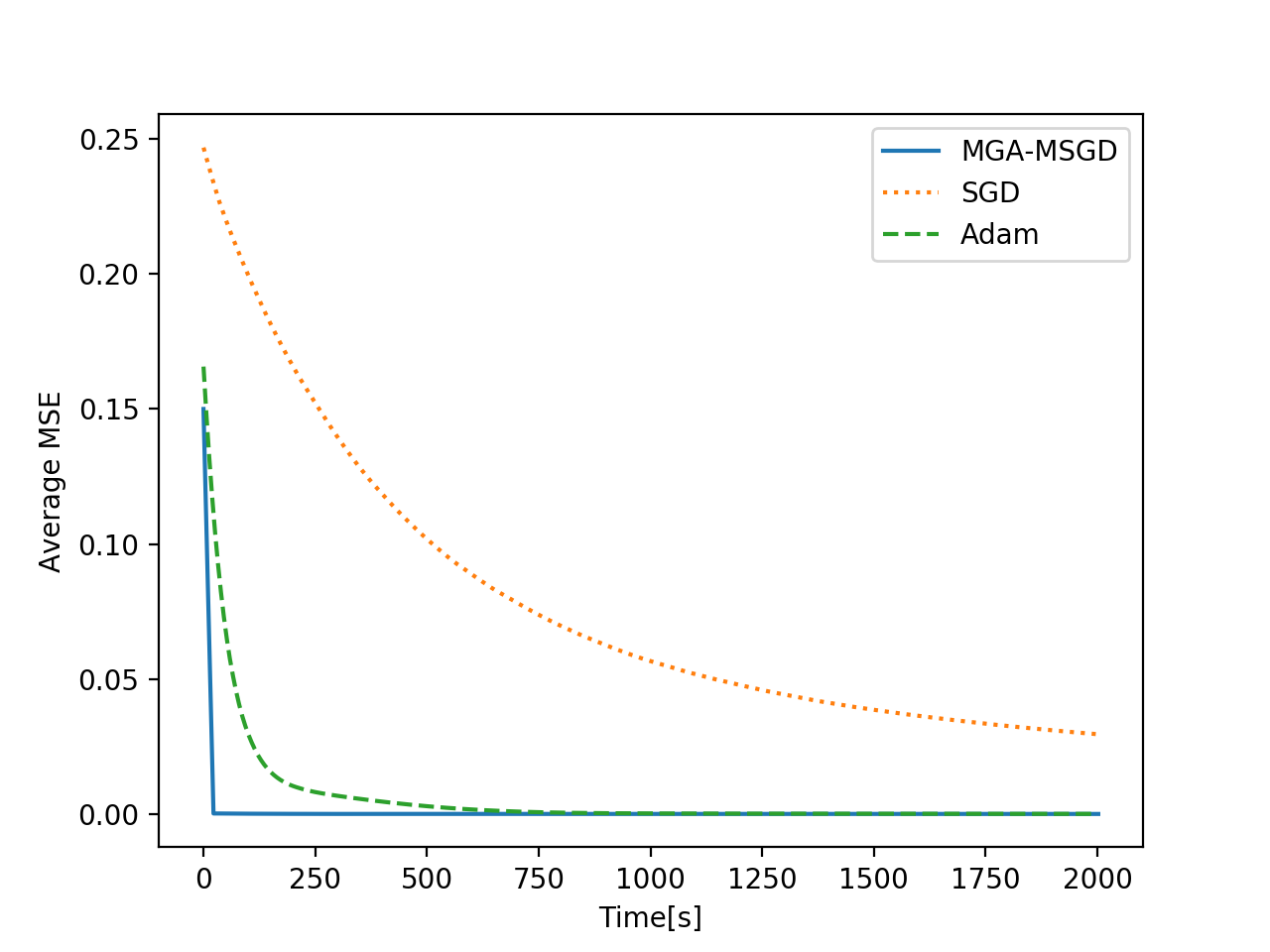}
	\centering
	\caption{The large scale picture of overall performance shows the efficiency of MGA-MSGD training algorithm and indicate that the major competitor is Adam optimiser}
	\label{fig_Performance}
	\end{subfigure} \quad 
	\begin{subfigure}{5.8cm}
	\includegraphics[width=5.8cm]{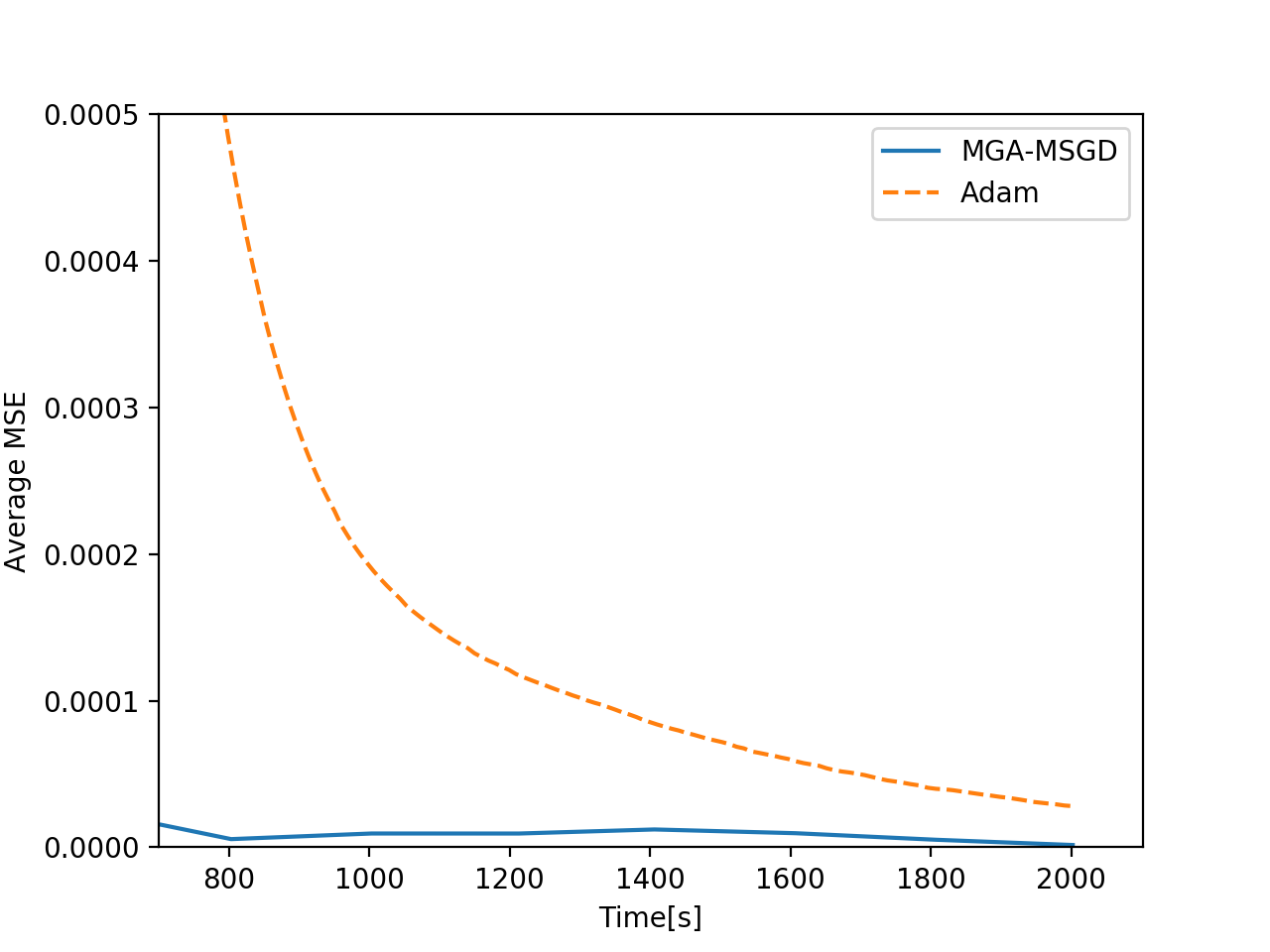}
	\centering
	\caption{The limited plot frame shows that even at very long times such as 1000[s] and 2000[s] MGA-MSGD outperforms the standard training approaches.}
	\label{fig_Performance_zoom}
	\end{subfigure}
	\caption{A performance comparison between the presented framework and two standard optimisers showing considerable improvement in both accuracy and efficiency.}
	\label{fig_Performance1}
\end{figure}

In the case of Adam optimiser the formulation is more complex where the parameters are updated via
\begin{equation}
\sbt^{(i)}_{kj}:=\sbt^{(i)}_{kj}- n \frac{\bar{m}}{\sqrt{\bar{v}}+ \bar{\epsilon}} \label{laststep} 
\end{equation}
where $\bar{\epsilon}$ is a term to improve the numerical stability. $\bar{m}$, $\bar{v}$ are defined and updated in each iteration via
\begin{equation}
\bar{m} := \frac{m}{1-\beta_1^t} \quad \quad \quad \quad \bar{v} := \frac{v}{1-\beta_2^t} 
\end{equation}
where, $t$ is the number of iteration and $\beta_1$ and $\beta_2$ are exponential decay rates for the moment estimates. The terms $m$ and $v$ are defined and updated in each iteration by
\begin{align}
m&:=\beta_1 m +(1-\beta_1) \frac{\partial \mathit{MSE}}{\partial \sbt^{(i)}_{kj}} \\ 
 v &:=\beta_2 v +(1-\beta_2) \left(\frac{\partial \mathit{MSE}}{\partial \sbt^{(i)}_{kj}}\right)^2 
\end{align}
We choose the initial values of $m$ and $v$ equal to zero.
Figure \ref{fig_Performance1} shows that MGA-MSGD converges to a minimum closer to the global minimum in remarkably smaller times. To our experience, the key to the higher efficiency is the use of a large learning rate which has become feasible by the help of the provided Modified Genetic Algorithm (MGA).

\subsection{Analysis of local displacement/output error}
\begin{figure}
\centering
	\begin{subfigure}{5.8cm}
	\includegraphics[width=5.8cm]{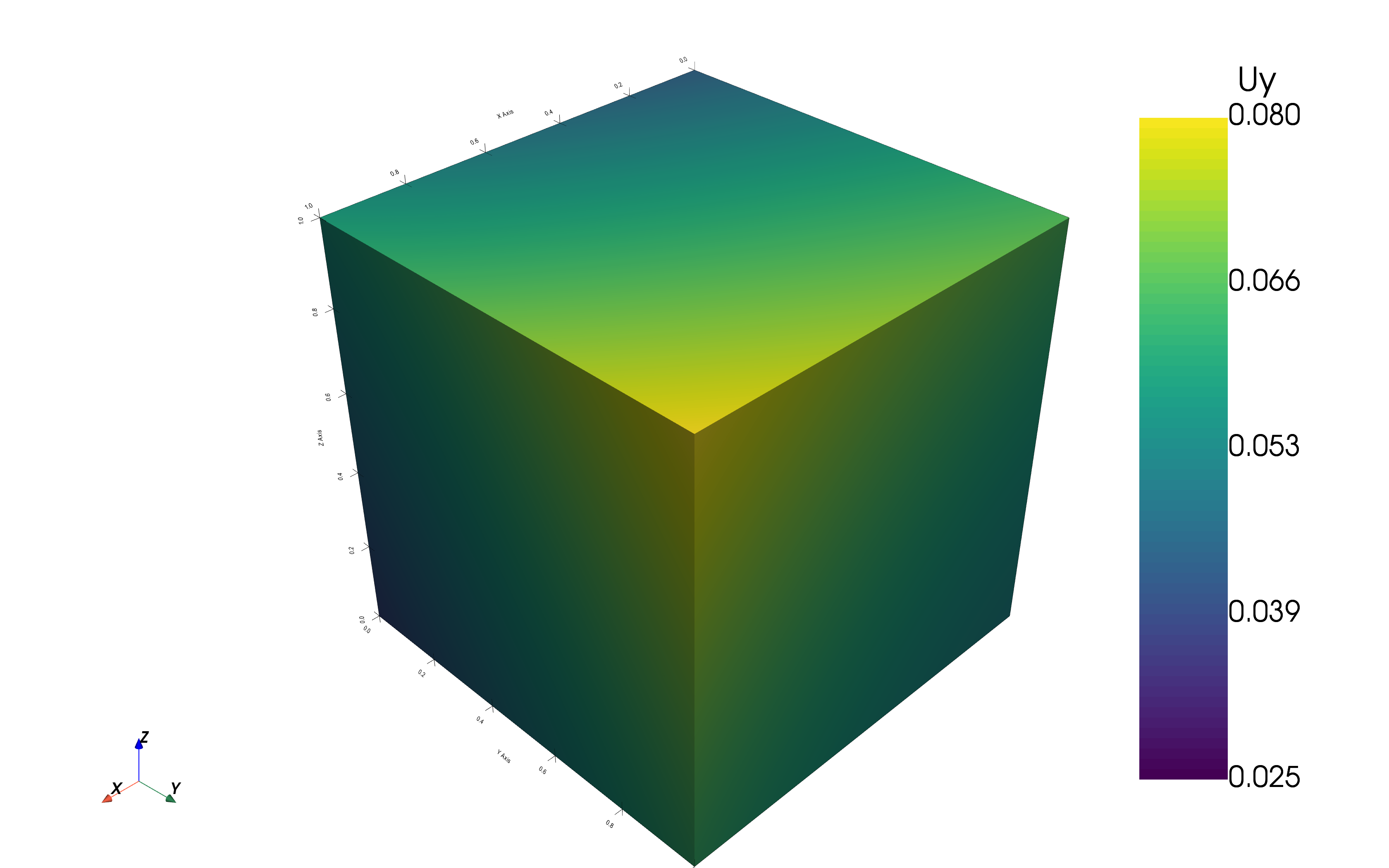}
	\centering
	\caption{}
	\label{fig_UyWrong}
	\end{subfigure} \quad 
	\begin{subfigure}{5.8cm}
	\includegraphics[width=5.8cm]{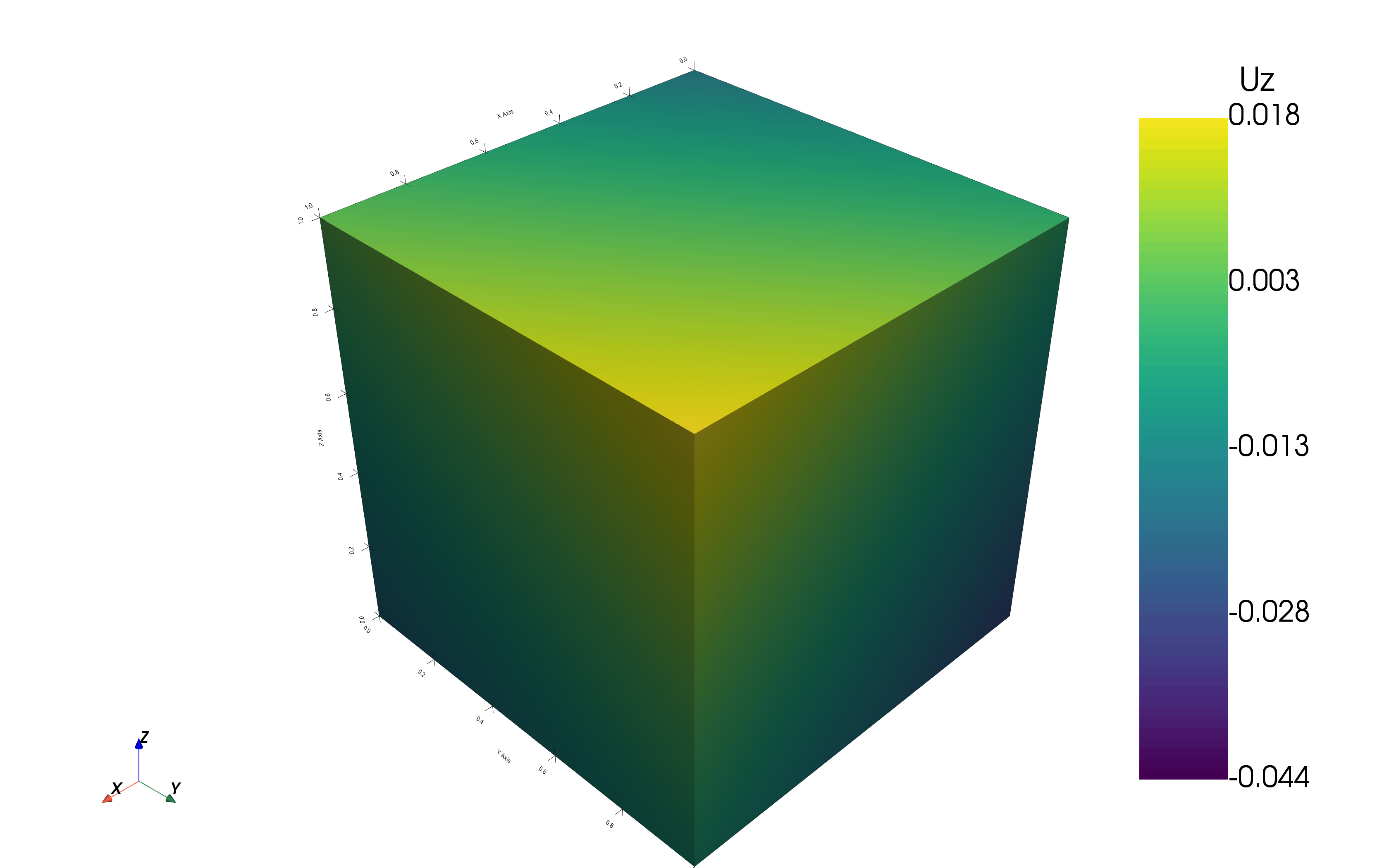}
	\centering
	\caption{}
	\label{fig_UzWrong}
	\end{subfigure}
	\caption{The contour plot of the resultant displacement components ($u_y$ and $u_x$) without enforcing Equations \eqref{normalisation} and \eqref{gamma}. We observe that the Balance of linear momentum in Equation \eqref{balance} is satisfied but the Dirichlet BC and the uniqueness condition in Equations \eqref{EqDBC}, \eqref{uniquex}, and \eqref{uniquey} are not imposed correctly. In other words, the exact results (provided in Figures \ref{fig_FEMUy} and \ref{fig_FEMUz}) can be obtained by a rigid body motion.}
	\label{fig_wrong}
\end{figure}

Apart from the statistical measures, it is important to see if the results make sense from a mechanical point of view as well. Figure \ref{fig_wrong} is the contour plot of the displacements in $\vett y$ direction of an analysis with acceptable $\mathit{MSE}$ loss, which, from a statistical viewpoint, is acceptable. However, from a mechanical viewpoint, to obtain the accurate response a rigid body motion is required which is also observed in similar works in the literature \cite{Raissi2018}.

\begin{figure}
\centering
	\begin{subfigure}{5.8cm}
	\includegraphics[width=5.8cm]{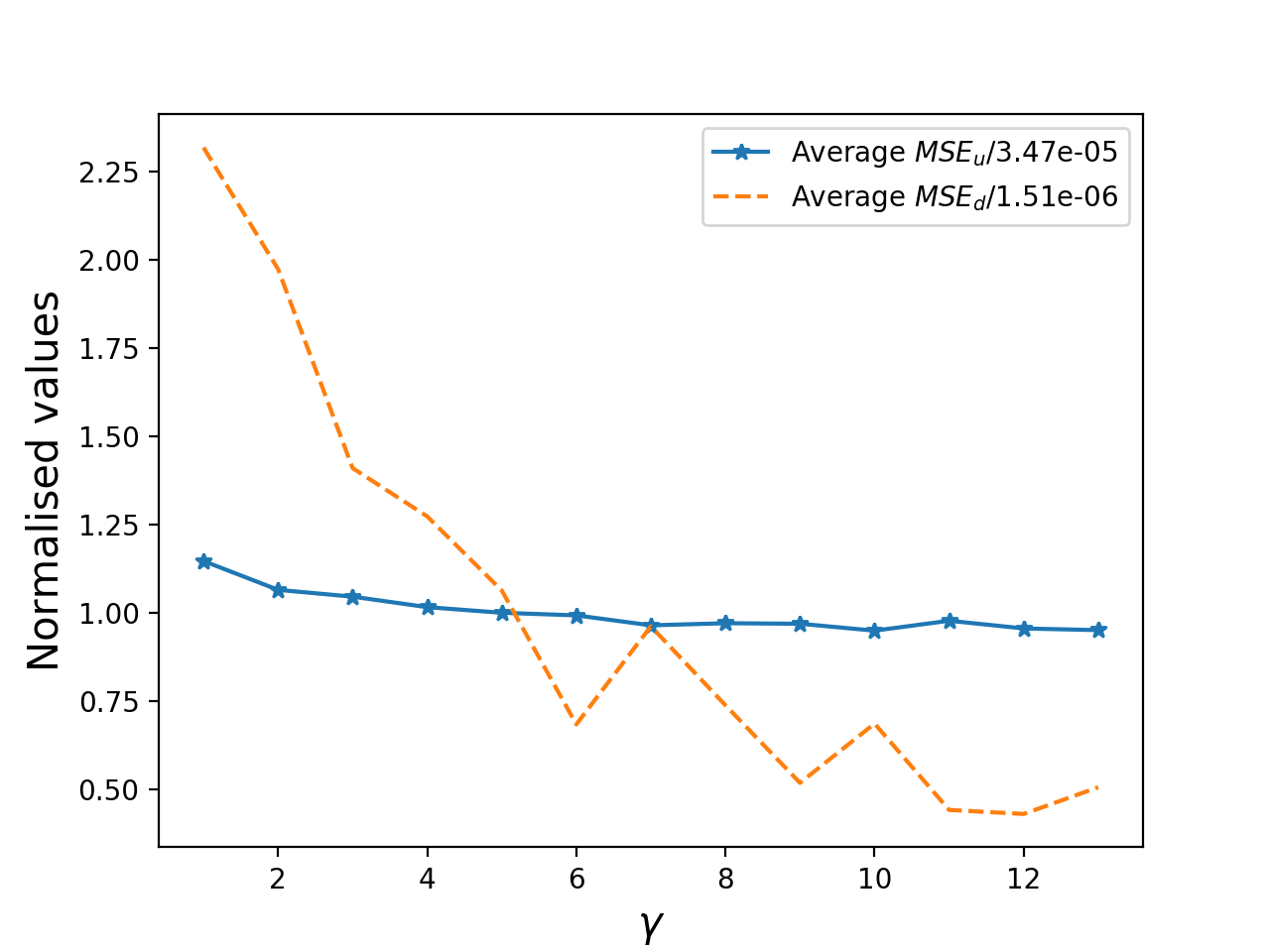}
	\centering
	\caption{An increase in $\gamma$ results in considerable decrease in $\mathit{MSE}_d$.}
	\vspace{0.3cm}
	\label{fig_gamma_MSEDisp_DBC}
	\end{subfigure} \quad 
	\begin{subfigure}{5.8cm}
	\includegraphics[width=5.8cm]{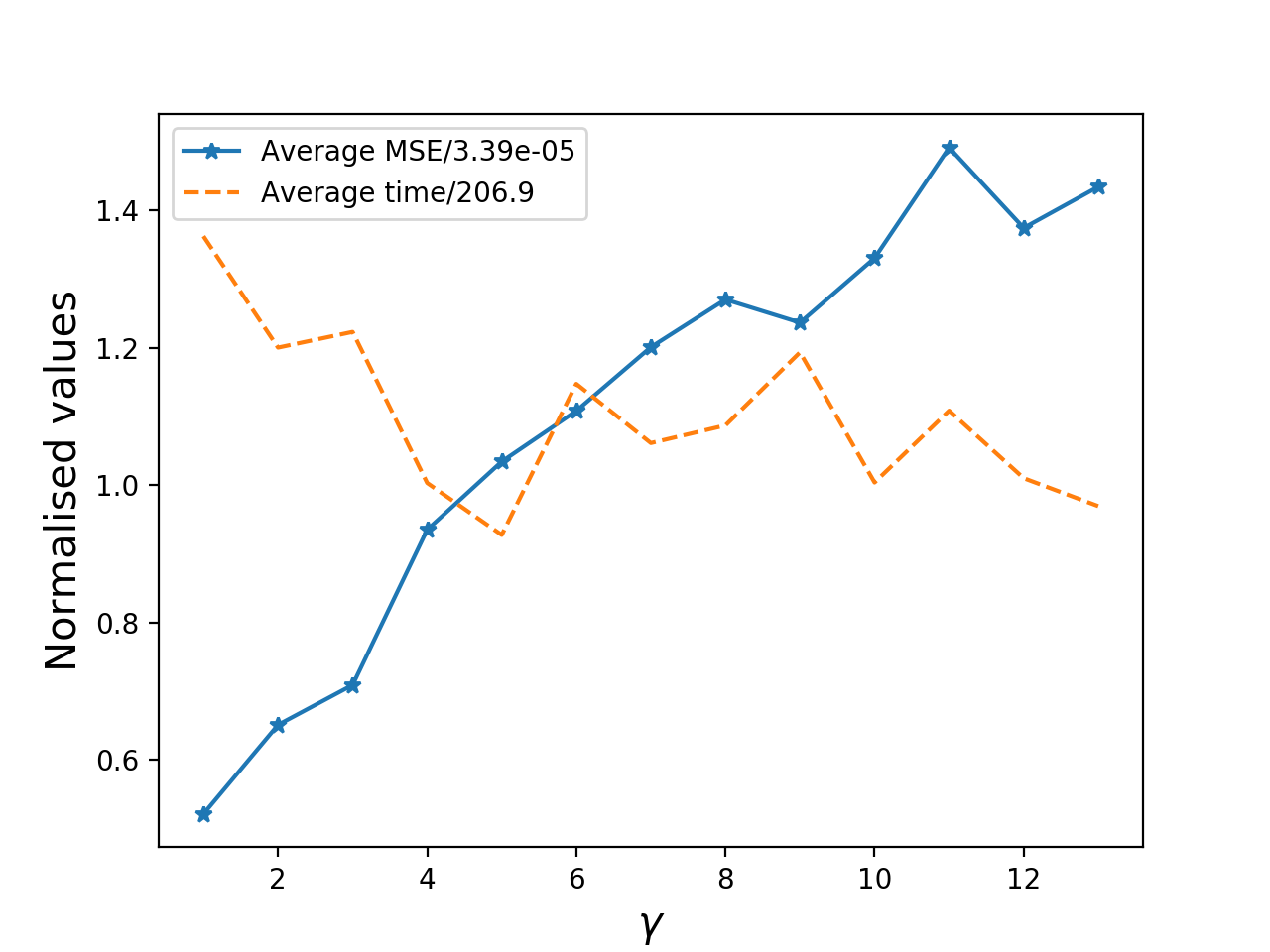}
	\centering
	\caption{The increase in $\mathit{MSE}$ is due to the multiplication of the distance function with the coefficient $\gamma$.}
	\label{fig_gamma_MSE}
	\end{subfigure}
	\caption{The parameter $\gamma$, which is a weight for strong enforcement of Boundary Conditions has a significant effect on the accuracy from a mechanical viewpoint, although it might cause an increase in the general $\mathit{MSE}$.}
	\label{fig_GMDB}
\end{figure}

In fact, although the balance of linear momentum seems to be enforced appropriately, the Dirichlet BC and uniqueness conditions are not represented correctly. The reason roots in the definition of $\mathit{MSE}$ (which is not a normalised/dimensionless value) as the values of stress components and their divergence can be considerably higher than the displacements, so the ANN training procedure only tends to satisfy the balance of linear momentum to a considerably higher extent than Dirichlet BC and uniqueness conditions. This problem can be solved  by normalising the stress $\vett{\sigma}$ as below
\begin{equation}
\hat{\vett{\sigma}} = \frac{\vett{\sigma}}{\lambda+2G}\label{normalisation}
\end{equation}
where, $\lambda$ and $G$ are the material parameters (Lam\'e constants).
Moreover, in order to ensure a strong/exact enforcement of Dirichlet BC we introduce the weight $\gamma$ as follows
\begin{equation}
\vett w_5(\vett x) = \vett r_5 (\vett x)*\gamma \quad  \quad \,on \quad  \quad \Gamma_d, \label{gamma}
\end{equation}
where $\gamma=0.05 n=6.25$ is chosen experimentally based on Figure \ref{fig_GMDB}. The use of the total number of points in the volume ($n$) in obtaining the parameter $\gamma$ is to ensure that the Dirichlet BCs are considered while it does not undermine other terms in the cost function $\mathit{MSE}$.
The mentioned conditions result in a more accurate mechanical response which is shown (of the same displacement components) in Figures \ref{fig_UYcomp} and  \ref{fig_UZcomp}.

At this stage, we solve the problem with the Finite Element Method (FEM) in order to obtain a reference solution that allows to study the local displacement error. This provides us with a criterion to study the resultant mechanical response. In fact, having the reference mechanical response of the problem ($\vett u$ in $\Omega$) we are able to quantify the distance function directly by
\begin{equation}
\mathit{MSE}_u =  \frac{1}{n}\sum_{i = 1}^{n} \vett w_7(\vett x_i) \cdot \vett r_7 (\vett x_i) \quad  \quad \,in \quad  \quad \Omega \label{eqMSEu}
\end{equation}
where 
\begin{align}
\vett r_7(\vett x) &= \vett u(\vett x) - \vett u_r (\vett x)  \quad  \quad \,in \quad  \quad \Omega\\
\vett w_7(\vett x) &= \vett r_7(\vett x)   \quad  \quad\quad  \,\,\quad  \quad \,in \quad  \quad \Omega
\end{align}
and where ${\vett{u}}_{r}(\vett x)$ is the reference results, to ensure that our results are accurate. Note that the latter is only used to test the results with no role in the network training procedure.

 \begin{figure}
\centering
	\begin{subfigure}{5.8cm}
	\includegraphics[width=5.8cm]{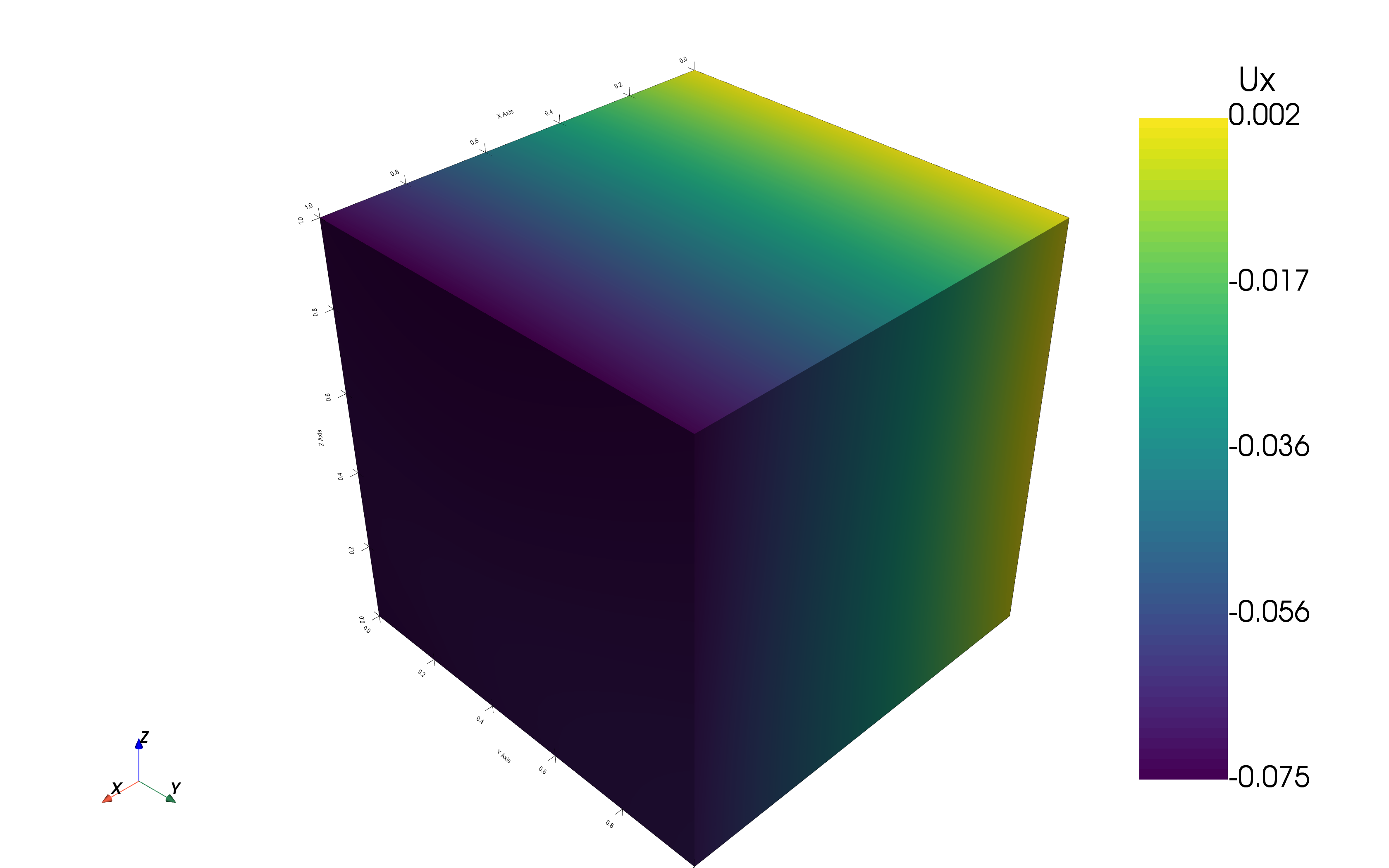}
	\centering
	\caption{$u_x$ from ANN.}
	\label{fig_Ux}
	\end{subfigure} \quad 
	\begin{subfigure}{5.8cm}
	\includegraphics[width=5.8cm]{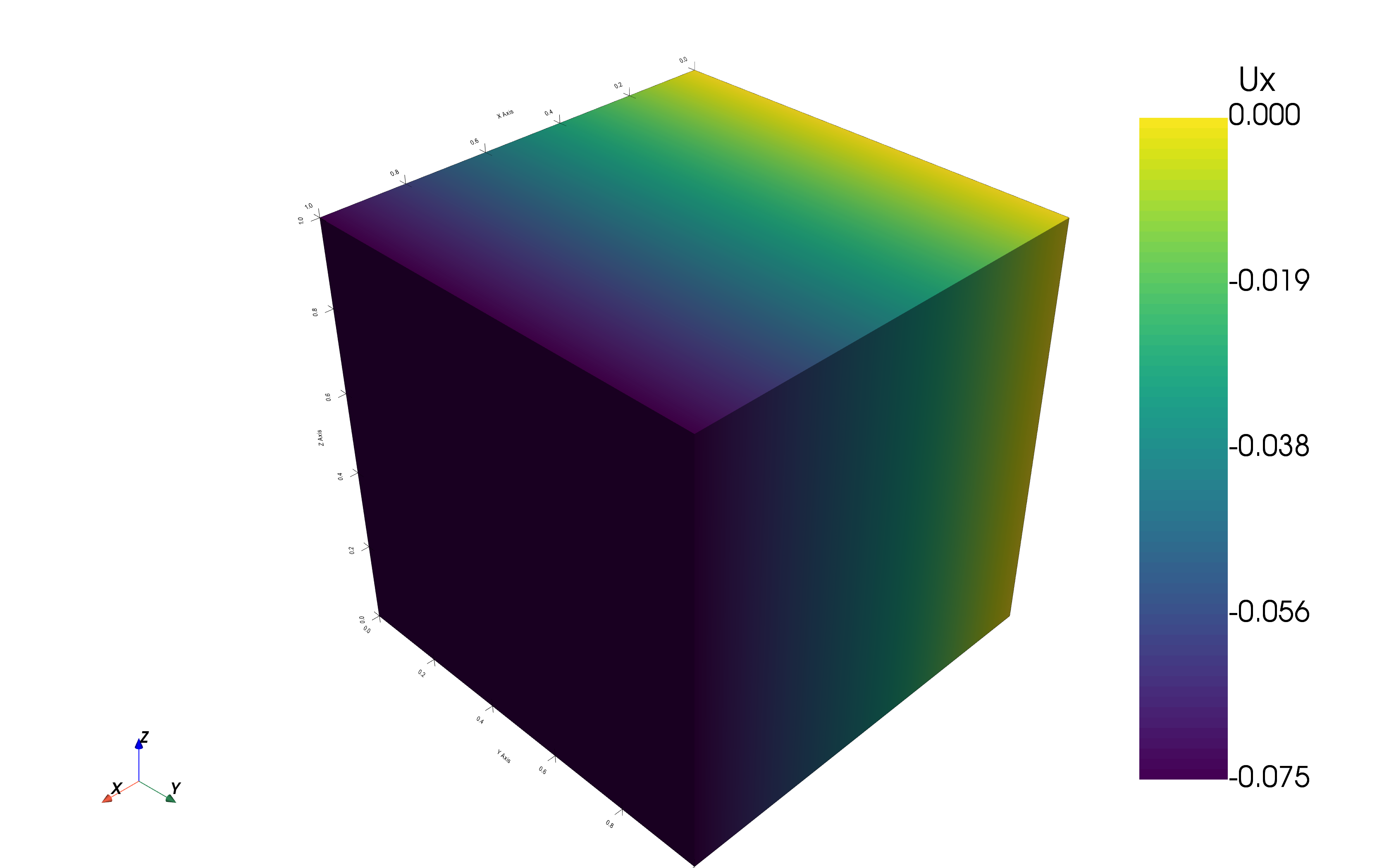}
	\centering
	\caption{Reference $u_x$.}
	\label{fig_FEMUx}
	\end{subfigure}
	\caption{The mechanical response provided by ANN trained via MGA-MSGD compared to the reference one provided by FEM.}
	\label{fig_UXcomp}
\end{figure}

 \begin{figure}
\centering
	\begin{subfigure}{5.8cm}
	\includegraphics[width=5.8cm]{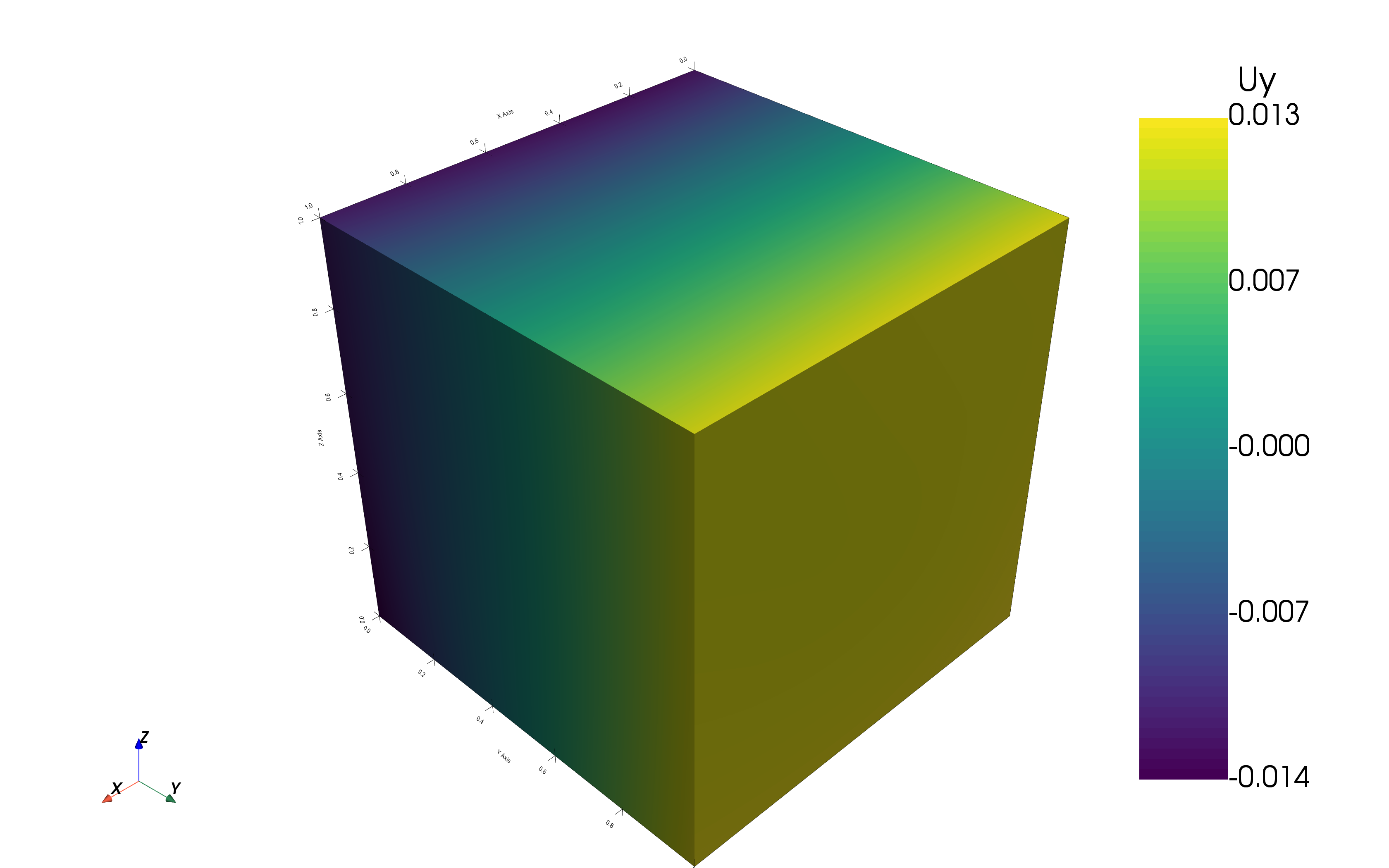}
	\centering
	\caption{$u_y$ from ANN}
	\label{fig_Uy}
	\end{subfigure} \quad 
	\begin{subfigure}{5.8cm}
	\includegraphics[width=5.8cm]{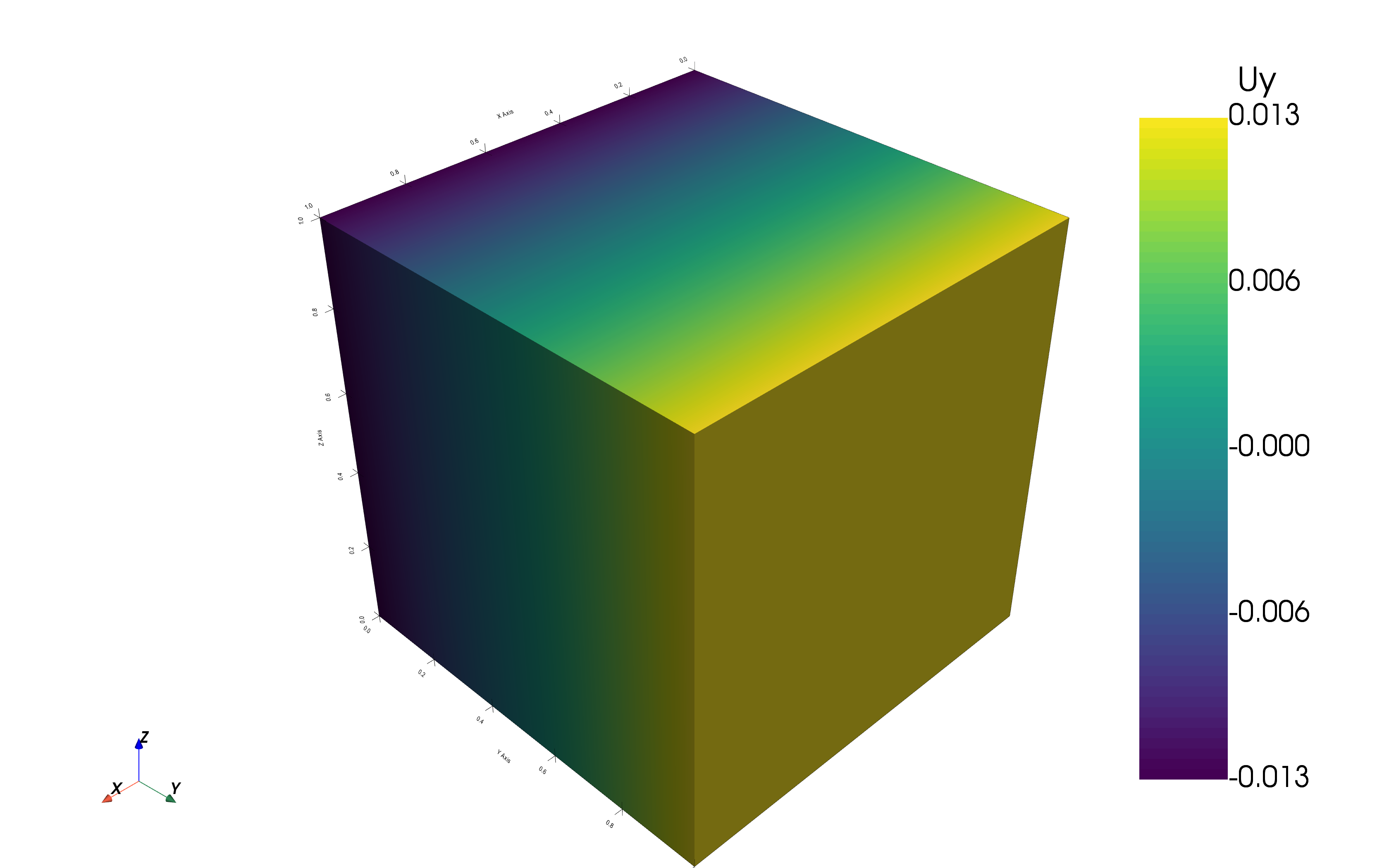}
	\centering
	\caption{Reference $u_y$.}
	\label{fig_FEMUy}
	\end{subfigure}
	\caption{The ANN results are in a good agreement with the reference ones showing the obtained accuracy.}
	\label{fig_UYcomp}
\end{figure}

 \begin{figure}
\centering
	\begin{subfigure}{5.8cm}
	\includegraphics[width=5.8cm]{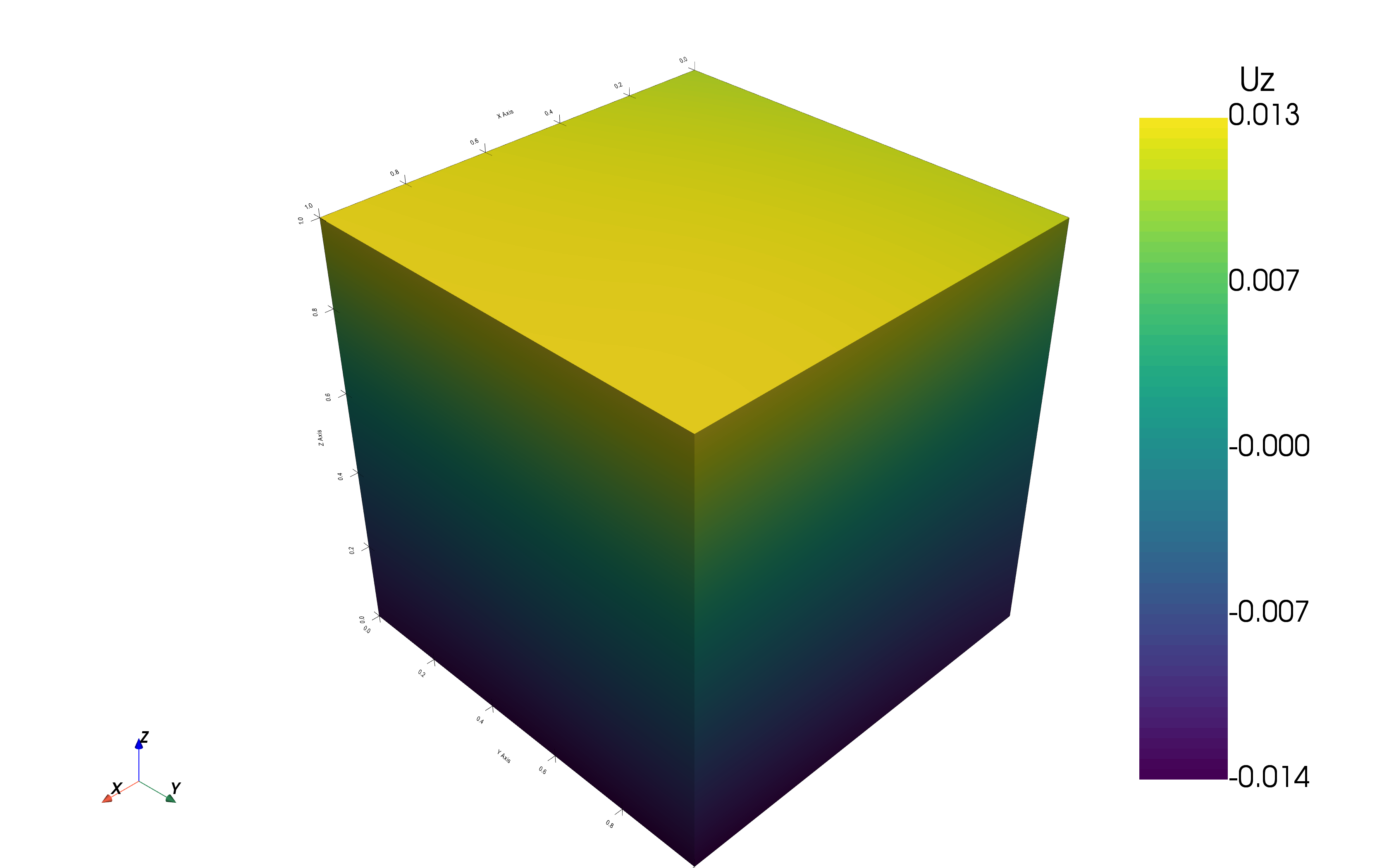}
	\centering
	\caption{$u_z$ from ANN}
	\label{fig_Uz}
	\end{subfigure} \quad 
	\begin{subfigure}{5.8cm}
	\includegraphics[width=5.8cm]{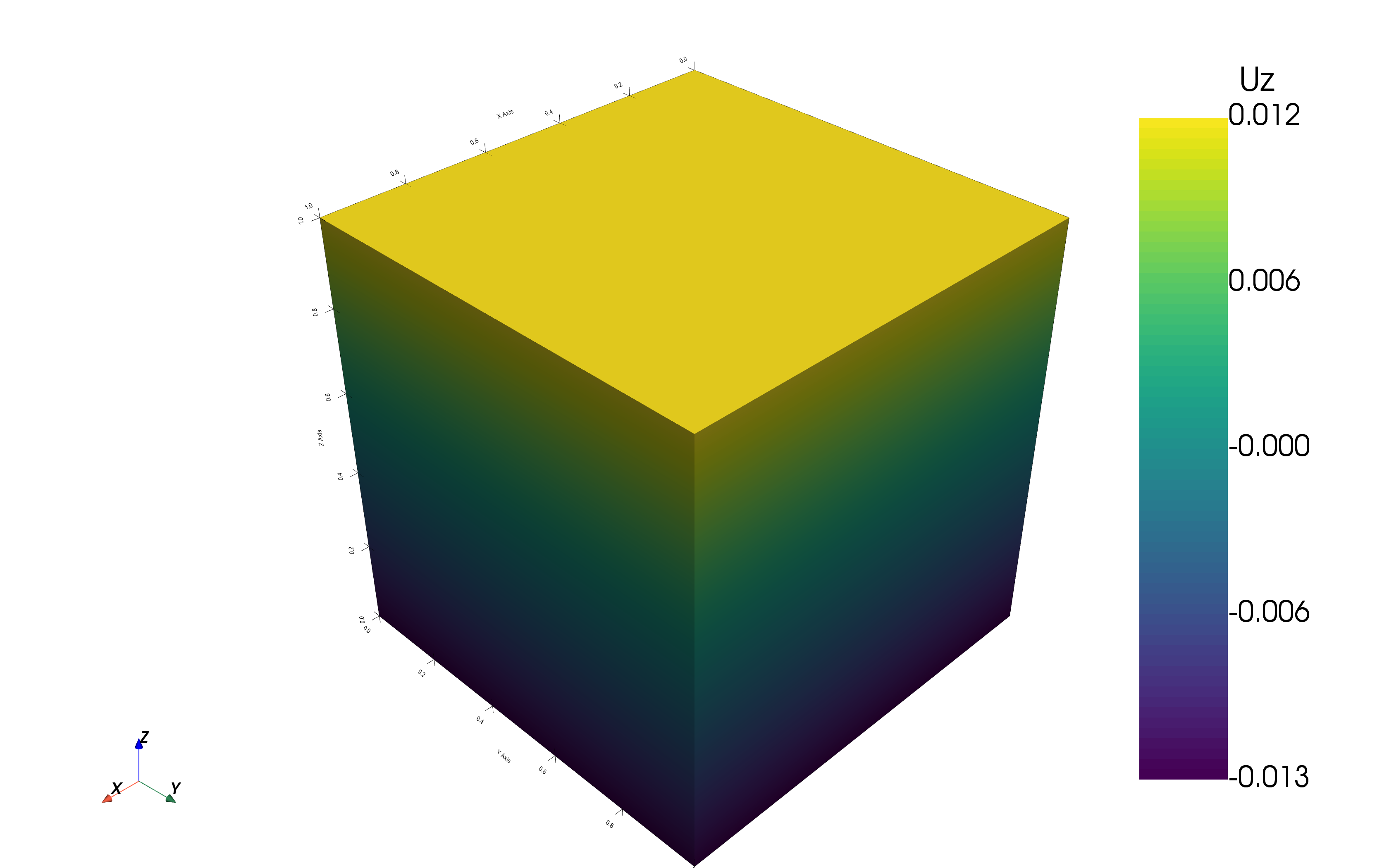}
	\centering
	\caption{Reference $u_z$.}
	\label{fig_FEMUz}
	\end{subfigure}
	\caption{Displacements in $\vett z$ direction from ANN and FEM as a mean of comparison.}
	\label{fig_UZcomp}
\end{figure}

A comparison between the solution provided by ANN and the reference displacements is shown in Figures \ref{fig_UXcomp}, \ref{fig_UYcomp}, and \ref{fig_UZcomp} shows that a good approximation of a 3D mechanical problem can be calculated by an ANN using the provided training strategy with an acceptable efficiency which was not possible via the conventional training algorithms.

In order to further evaluate the potential of the presented framework we solve the same boundary value problem with $\vett u_0 = (0,0,0)$ on $\Gamma_d$ which results in slightly more complex deformation. In this case, the displacements in $\vett y$ and $\vett z$ directions decrease approaching $\Gamma_d$.

 \begin{figure}
\centering
	\begin{subfigure}{5.8cm}
	\includegraphics[width=5.8cm]{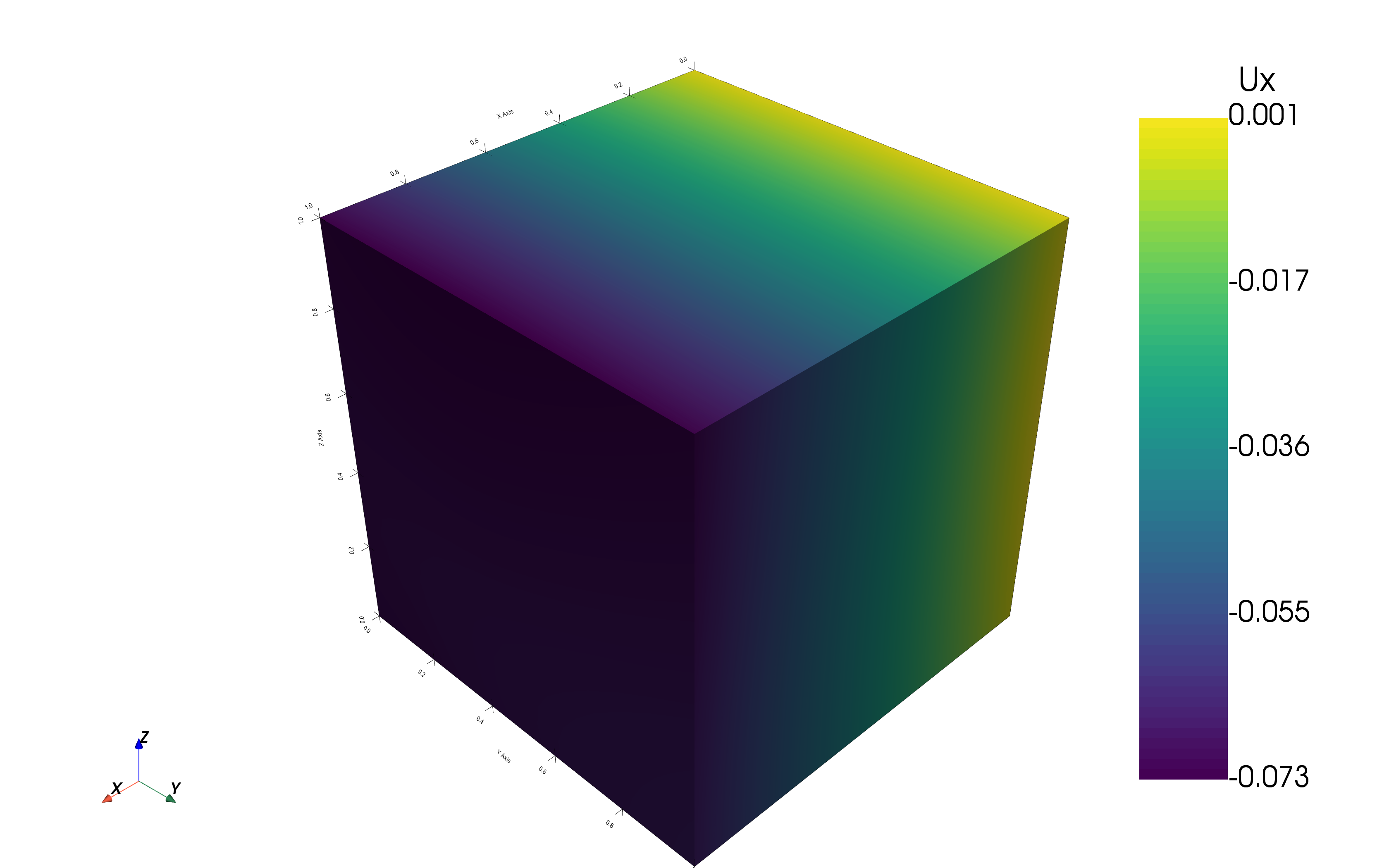}
	\centering
	\caption{$u_x$ from ANN.}
	\label{fig_Ux}
	\end{subfigure} \quad 
	\begin{subfigure}{5.8cm}
	\includegraphics[width=5.8cm]{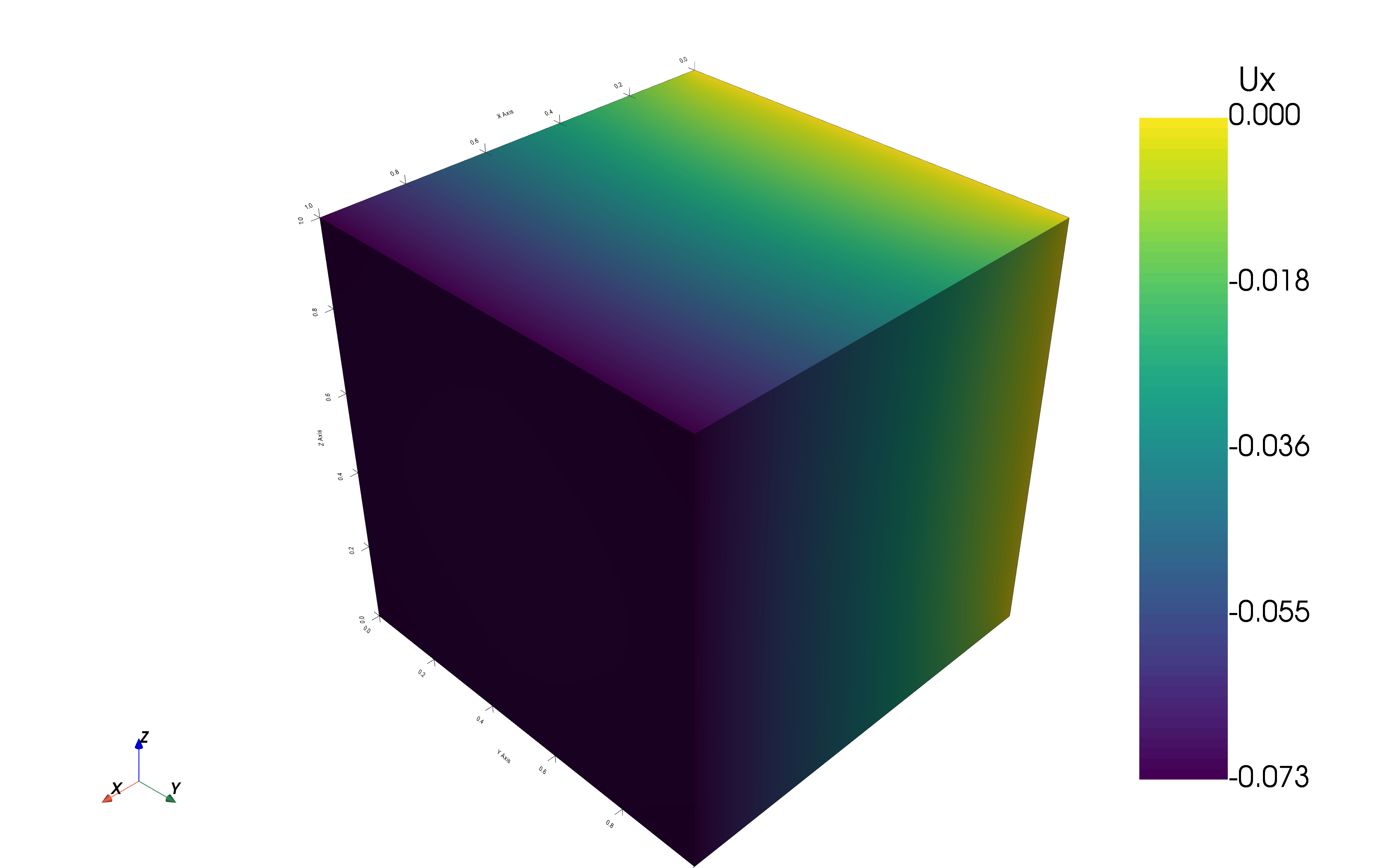}
	\centering
	\caption{Reference $u_x$.}
	\label{fig_FEMUx}
	\end{subfigure}
	\caption{
	The mechanical response imposing Dirichlet BC $\vett u_0 = (0,0,0)$ on $\Gamma_d$ (a) provided by ANN trained via MGA-MSGD (b) the reference one provided by FEM.}
	\label{fig_UXcomp}
\end{figure}

 \begin{figure}
\centering
	\begin{subfigure}{5.8cm}
	\includegraphics[width=5.8cm]{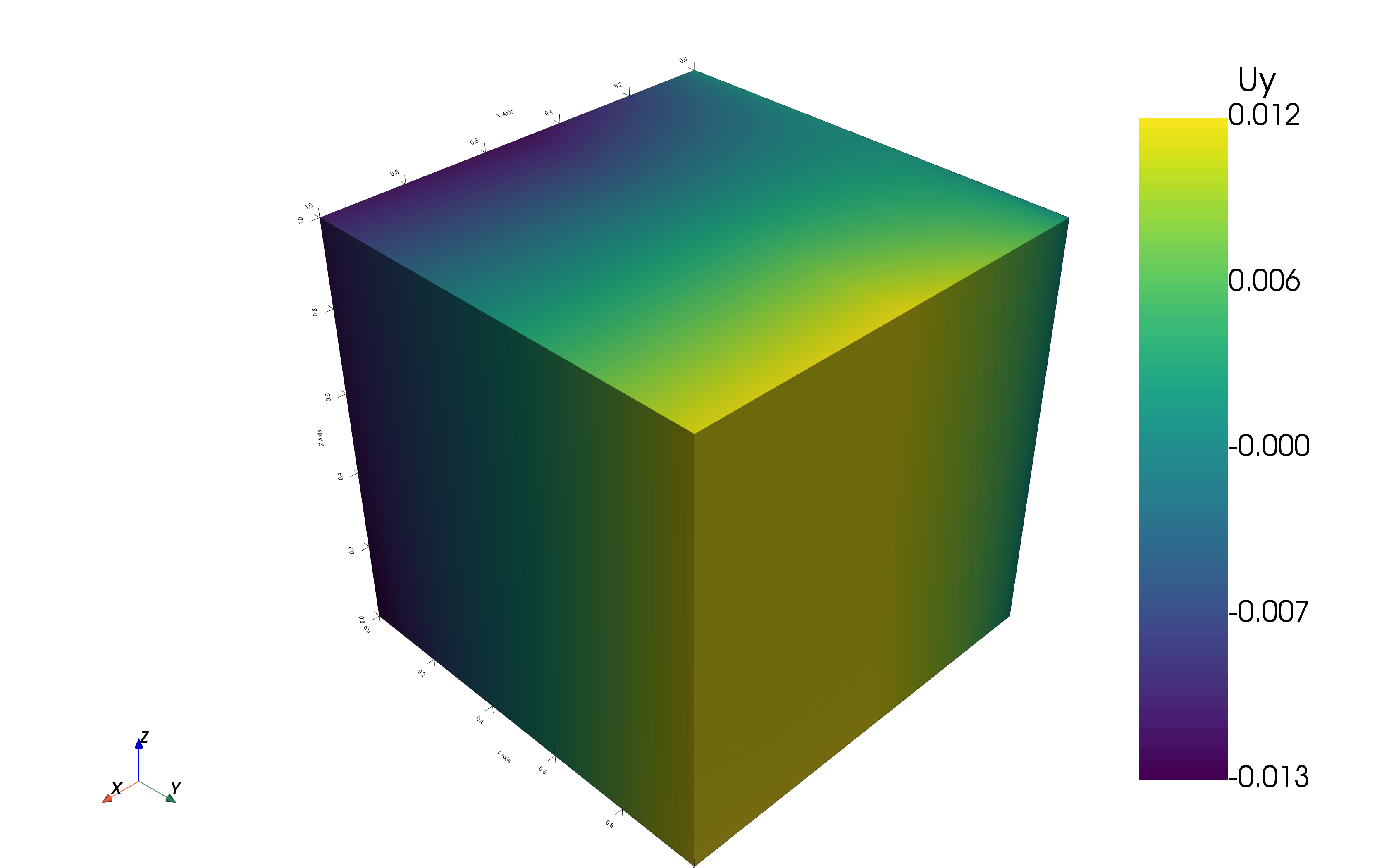}
	\centering
	\caption{$u_y$ from ANN}
	\label{fig_Uy}
	\end{subfigure} \quad 
	\begin{subfigure}{5.8cm}
	\includegraphics[width=5.8cm]{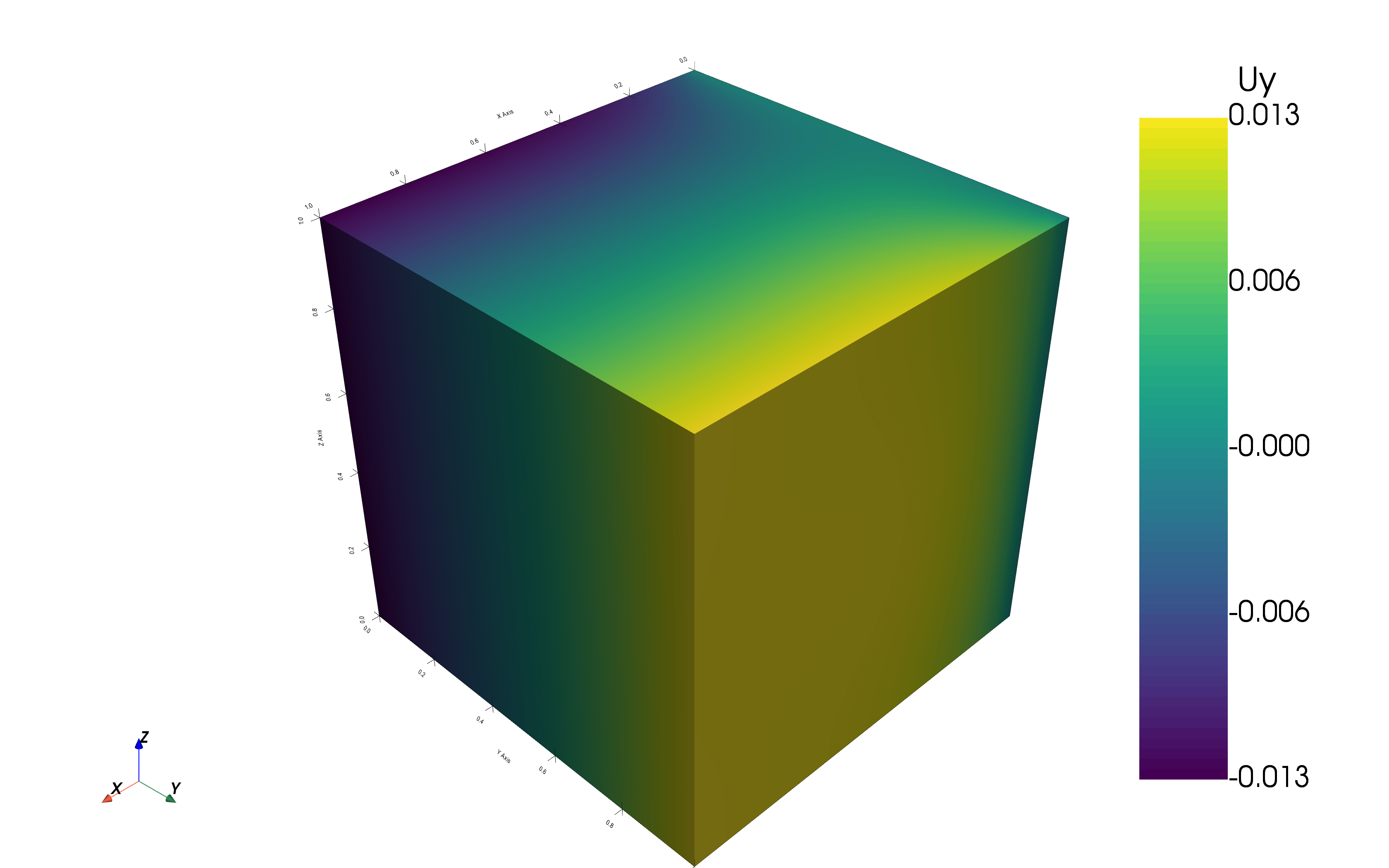}
	\centering
	\caption{Reference $u_y$.}
	\label{fig_FEMUy}
	\end{subfigure}
	\caption{The ANN results are in a good agreement with the reference ones highlighting its feasibility.}
	\label{fig_UYcomp}
\end{figure}

 \begin{figure}
\centering
	\begin{subfigure}{5.8cm}
	\includegraphics[width=5.8cm]{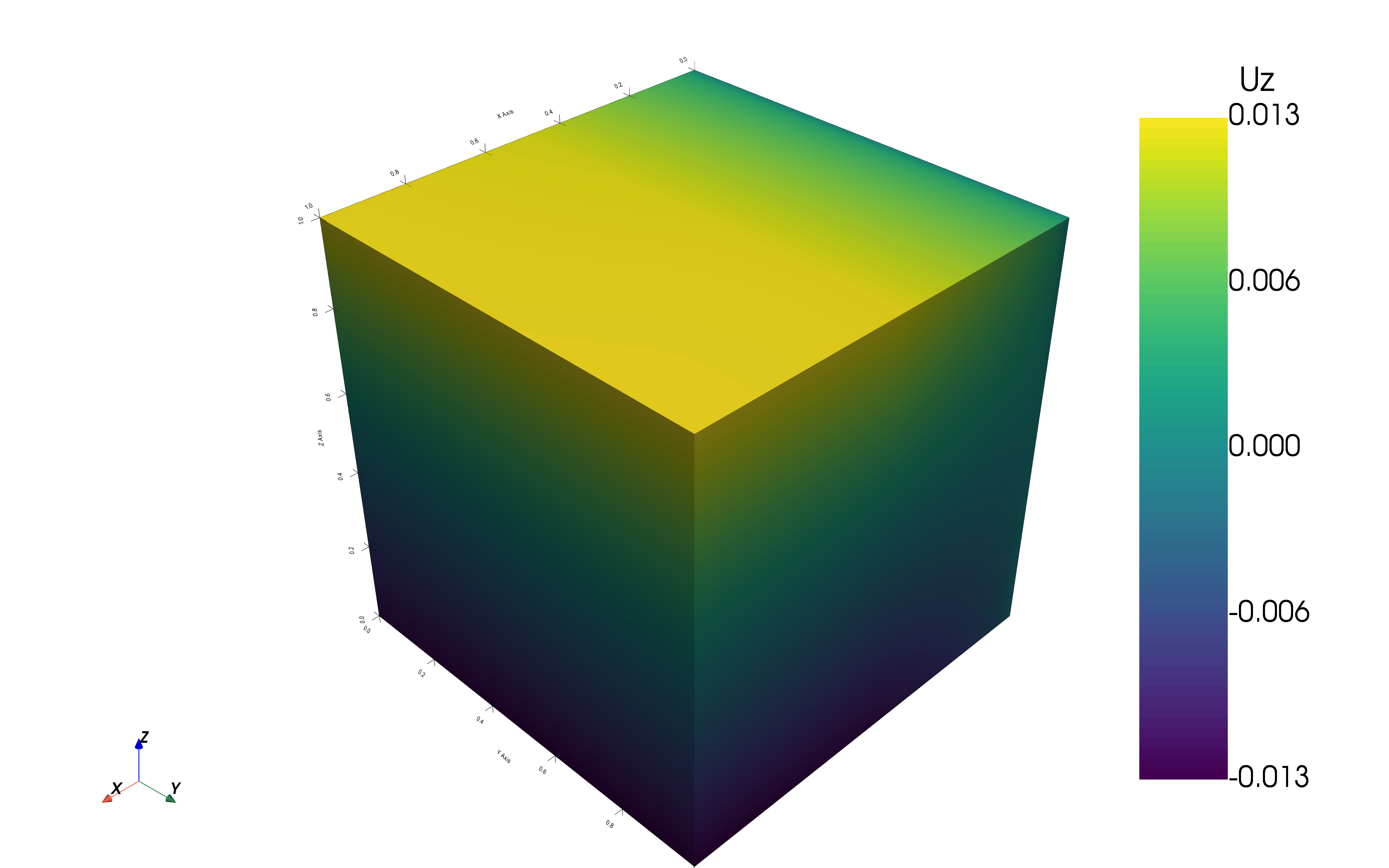}
	\centering
	\caption{$u_z$ from ANN}
	\label{fig_Uz}
	\end{subfigure} \quad 
	\begin{subfigure}{5.8cm}
	\includegraphics[width=5.8cm]{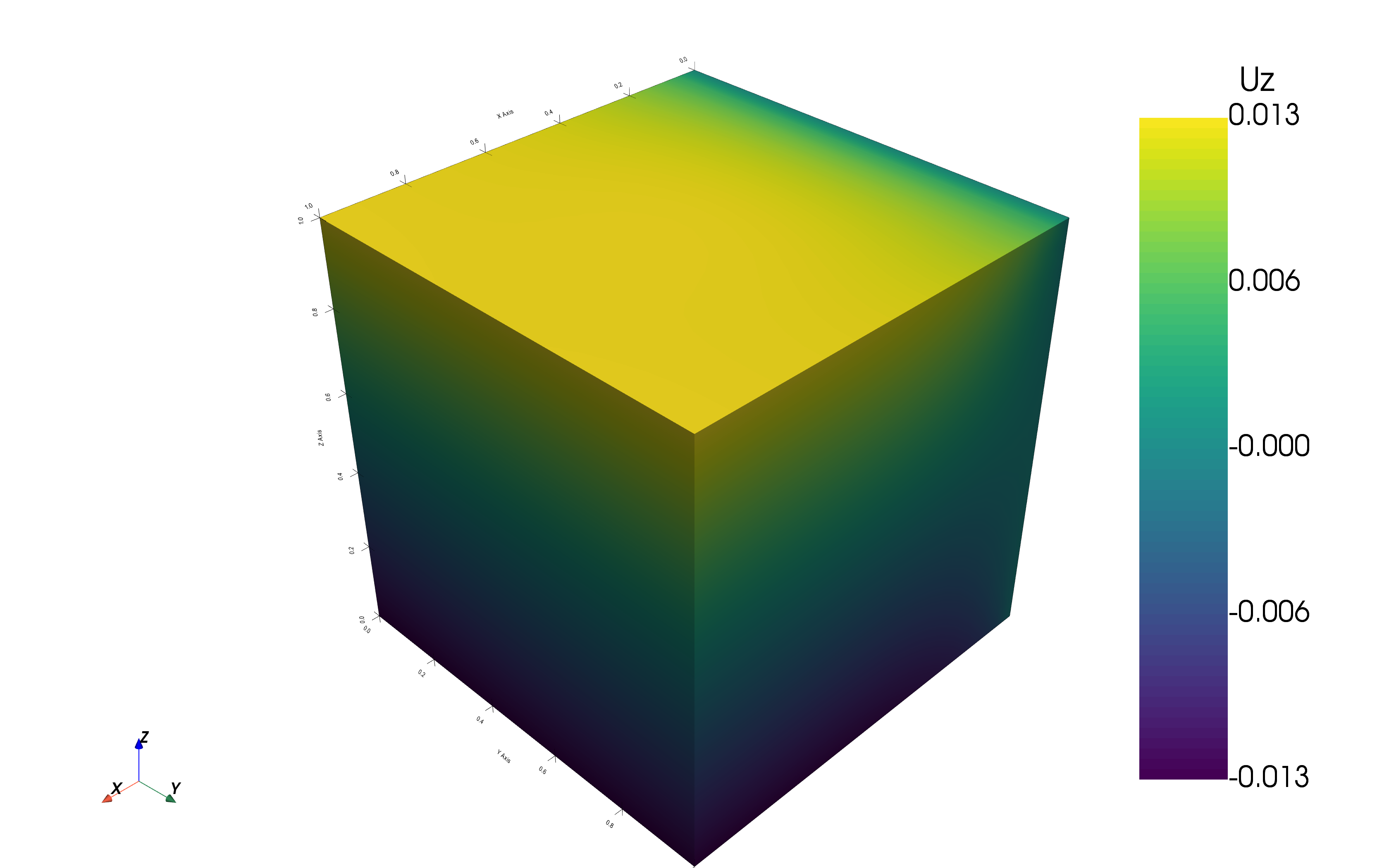}
	\centering
	\caption{Reference $u_z$.}
	\label{fig_FEMUz}
	\end{subfigure}
	\caption{Displacements in $\vett z$ direction are also in agreement with the reference one.}
	\label{fig_UZcomp}
\end{figure}

  \begin{figure}
  \centering
  \includegraphics[width = 8cm]{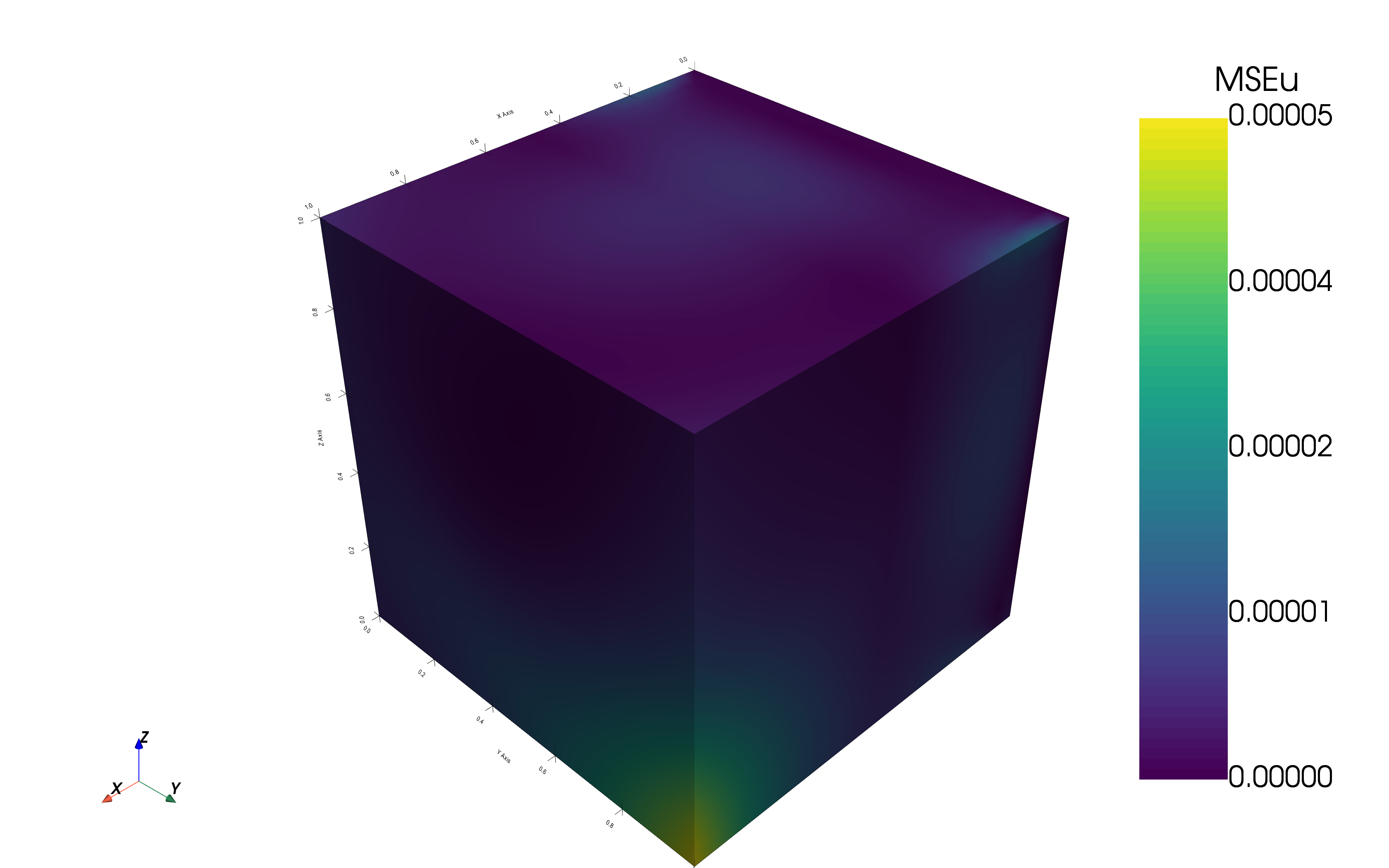}
  \caption{Local displacement error $\mathit{MSE_u}$ calculated via Equation \eqref{eqMSEu}}
  \label{fig_MSEu}
  \end{figure}
  
Finally, we visualise the local displacement error $\mathit{MSE_u}$ over the domain in Figure \ref{fig_MSEu} which shows that although the error is acceptable it exists mostly on the corners and where the displacement gradients are higher leaving the possibility to further refinements of the method in future works.

\section{Conclusion} \label{sec_conclusion}
We have introduced the hybrid MGA-MSGD training approach for ANNs with application to solving mechanical problems in three spatial dimensions described by elliptic PDEs. The latter application requires third order derivatives of the outputs with respect to the inputs and each network parameter, which is considerably time-consuming, especially, if we need a dense training data set/points and a large number of network's parameters.
This novel training approach includes modified GA adjusting the learnable parameters' components which cause the error explosion so that we can employ large learning rates (CSGD) avoiding several local minima. This is followed by a Fine-scale SGD training procedure to obtain the most accurate results. The method introduces some new non-mechanical parameters to be tuned via experiments and sensitivity analysis for a timely procedure. We show that the obtained results are less sensitive to the number and distribution of data points, and learning rates which are crucial achievements. Furthermore, this method allows us to exploit a larger potential of the network hence obtaining accurate results from small networks. The effectiveness of the training procedure is compared with two competitors, namely the classic SGD and another SGD-based optimiser (Adam optimiser), showing a significant improvement in both accuracy and efficiency. 
 The obtained results can be accurate up to a rigid-body motion. We have introduced the stress normalisation and displacement boundary condition enforcement weight in order to eliminate the need for any post-processing such as the mentioned rigid-body motion.
 The same problem is then solved via FEM which, having sufficiently fine space discretisation, is considered as the reference solution. The accurate final mechanical response obtained via the presented training method highlights its reliability. This approach considers the strong form of the governing equations providing results that are directly differentiable and free of space discretisation (mesh-free). In practice, in order to obtain the response of a specific point, one does not need to discretise the whole domain and obtain the response of all the nodes (path-independent response). This advantage reduces the complexity of the problem implementation while it provides us with more flexibility. The authors highlight that this method should not be taken as a replacement for approaches such as FEM and finite difference method but complementing them.
The potential of the presented framework is further studied by solving a slightly more complex problem highlighting that it can be applied to a broad range of the scenarios of interest such as poroelasticity  \citep{DEHGHANI2020103996, HdehghaniThesis, Hdehghani}, by including time as the 4th dimension, and non-linear elasticity with residual stress \citep{DEHGHANI201951, Njfiber, JHA201979} by using the corresponding governing equations and measures.

\section*{Acknowledgements}
We acknowledge the support of this research work via the framework of DTU DRIVEN, funded by the Luxembourg National Research Fund (PRIDE17/12252781), and the project CDE-HUB, funded by the Luxembourg Ministry of Economy (FEDER 2018-04-024).



\bibliographystyle{plain}
\bibliography{Dbib}

\end{document}